\documentclass[preprint, a4paper, 10pt]{elsarticle}

\usepackage{lmodern}
\usepackage{lineno,hyperref}
\usepackage[english]{babel}
\usepackage{amsmath,amsfonts,amssymb,amsthm}
\usepackage{graphics}
\newtheorem{Rem}{\large \sf Remark}[section]
\newtheorem{Lem}[Rem]{\large \sf Lemma}

\newtheorem{Prop}[Rem]{\large \sf Proposition}

\newtheorem{Cor}[Rem]{\large \sf Corollary}

\newtheorem{Convention}[Rem]{\large \sf Convention}

\newtheorem{Def}{\large \sf Definition}
\newtheorem*{ContToPL_lem}{Lemma \ref{ContinuousToPL_lem}}

\newcommand{\PR}{{\vspace*{-0.2cm}\sf{\em{Proof.}} } }
\newcommand{\EPR}{$\diamondsuit$\vspace*{0.3cm}}

\newcommand{\smbf} {\bf \small}

\begin{document}

\begin{frontmatter}
\bibliographystyle{elsarticle-num}

\title{Topological Laminations on Surfaces}

\author{Luis-Miguel Lopez}
\ead{lopez@ed.tokyo-fukushi.ac.jp}

\address{Tokyo University of Social Welfare, Dept. of Education, \\
Gunma Pref., Isesaki City, Sannoncho 2020-1, 372-0831 Japan}

\begin{abstract} We present and prove a topological characterization of geodesic 
laminations on hyperbolic surfaces of finite type.
\end{abstract}


\end{frontmatter}

\section{Introduction}
{\em Geodesic laminations} on hyperbolic surfaces were defined by Thurston in the 
early seventies~\cite[Sect. 8.5]{Thu80}, and have proved since to be an 
important tool in low-dimensional topology~\cite{Thu80,CB88,Bon01,BJCRMY12}: 

\begin{Def}
A {\smbf geodesic lamination} $\cal L$ on a hyperbolic surface $\Sigma$ is a
non-empty closed subset of~$\Sigma$ which is a union of simple and pairwise
disjoint complete geodesics.
\label{geodlam_def}
\end{Def}

\noindent
This definition is short and simple. However, the geodesics in a lamination, 
seen as curves, have the following important homotopical property: they can be 
continuously deformed onto embedded finite graphs, satisfying specific 
differentiability conditions and called 
{\em train tracks}, that is, geodesic laminations are {\em carried by} 
train-tracks~\cite{Thu80, PH93}. 
This property led Hatcher to redefine a measured lamination as 
a curve set carried by a train track taken among finitely many predefined 
ones, in order to directly obtain results independent of the choice of a 
metric~\cite{Hat88}. Other metric independent results have been 
proved more recently~\cite{ZB04}, again using this property. Behind these 
results is the fact that geodesic laminations can be topologically 
characterized, a characterization that we make explicit as follows: 

\begin{Def}
A {\smbf topological lamination} $\cal T$ on a surface $\Sigma$ is a non-empty 
subset of $\Sigma$ which is a union of simple and pairwise disjoint curves, all 
closed or two-way infinite, such that the curves are pairwise non-homotopic, 
they are all carried by a common train track $\Gamma$, and $\cal T$ is maximal 
for this~$\Gamma$ and remains so after any isotopies applied to the curves. 
\label{lam_def}
\end{Def}

\noindent
The above definition is about ordinary continuous curves, and the 
train track condition is combinatorial, as train tracks can be seen as plain embedded 
graphs together with vertex crossing rules for the carried paths. Moreover, quite 
remarkably, the topological hypothesis of being closed in Definition~\ref{geodlam_def} 
becomes the setwise hypothesis of being maximal with respect to set inclusion 
in Definition~\ref{lam_def}. 

Now, although it is admitted that geodesic laminations can be characterized up to 
isotopy as in Definition~\ref{lam_def} (see~\cite[Prop. 8.9.4]{Thu80} for a  
result in this direction), to the knowledge of the author this fact has 
not been proved, not even explicitly written down. 
Such a characterization shows nevertheless some advantages. 
Not only maximality gives an alternative 
way of checking that a set of curves is a lamination, but also a combinatorial 
viewpoint of laminations is made visible and can be exploited for itself. 
In particular, by labeling the edges of the train tracks, laminations can be made 
into usual {\em shifts}, creating bridges 
towards symbolic dynamics and word combinatorics~\cite{FZ08, LN12, N14, BDDPRR16}. 

The aim of this paper is to present a detailed exposition of the proof of 
the equivalence up to isotopy between Definitions~\ref{geodlam_def} and~\ref{lam_def}. 
Noteworthily, in the course of this proof two issues specific to topological laminations 
are raised: Defining a notion of {\em complementary region}, and dealing 
with non-compact leaves homotopic to compact curves. 
We also prove that a topological lamination carried by a train track 
$\Gamma$ can be moved by an ambient isotopy into a regular neighborhood of $\Gamma$ 
(see Proposition~\ref{Engulfing_prop}). 
This latter result says that laminations can be described with respect to fat graphs 
(equiv. ribbon graphs) \cite{LZ04, EM13}, relaxing their dependency to the surfaces 
themselves, and thus making laminations even more combinatorial objects. 

\section{Basic Definitions}
\label{lam_intro}
\subsection{Surfaces and PL structures}
A {\smbf surface of finite type} is a closed surface from which 
finitely many points, the {\smbf punctures}, have been removed. 
All the surfaces of finite type except the torus and the sphere 
admit complete metrics with constant curvature~$-1$, making them 
{\smbf hyperbolic surfaces}, but only those whose area is 
finite are called {\smbf hyperbolic surfaces of finite type}. 
Hyperbolic surfaces all share the same universal 
covering space, the {\smbf hyperbolic plane} $\mathbb{H}^2$, so 
they are quotients of $\mathbb{H}^2$ by discrete 
subgroups of isometries acting properly and discontinuously, without 
fixed points in $\mathbb{H}^2$. A subsurface, closed as a subset of a 
complete hyperbolic surface and whose frontier components are  
geodesics, is called a {\smbf hyperbolic surface with geodesic 
boundary}, and so is any surface isometric (in the Riemannian sense) 
to it. 

A topological manifold is {\em triangulable} if it is homeomorphic to 
a simplicial complex. When in addition the simplicial complex 
admits an atlas such that one can pass from chart to chart by piecewise 
linear functions, the manifold together with the complex is called a 
{\smbf PL manifold}, and one says the triangulation {\smbf is PL}. A 
piecewise linear map between PL manifolds is called a {\smbf PL map}. 
It is known that every manifold of dimension no more than~3 is triangulable 
and that every triangulation is PL in a unique way up to PL homeomorphism. 
Thus a tessellation of a hyperbolic surface of finite type $\Sigma$ 
by hyperbolic triangles defines a PL structure, which is unique, 
though there are many non-equivalent hyperbolic structures on $\Sigma$. 
Likewise, the universal covering by $\mathbb{H}^2$ of $\Sigma$, gives it a 
{\em differentiable}, or {\em smooth} structure which is unique up to 
diffeomorphism. In the sequel $\Sigma$ is assumed to be endowed with both 
its PL and smooth structures. 

\subsection{Curves, Geodesics and Laminations}
\label{CGL_sub}
A ({\smbf PL}) {\smbf curve}~$\gamma$ on a surface $\Sigma$ is a (PL) continuous map, 
either from a closed connected subset $J$ of the real line, or from the circle to 
$\Sigma$. In the latter case $\gamma$ is said to be {\smbf closed}; if the map defining 
the curve is injective, $\gamma$ is said to be {\smbf simple}; if $\gamma$ is simple, 
closed, and bounds neither a disk nor a once punctured disk it is said to be 
{\smbf essential}; if~$J=\mathbb{R}$, $\gamma$ is said to be {\smbf two-way infinite}; 
if~$J=\mathbb{R}_{\ge 0}$ or~$J=\mathbb{R}_{\le 0}$ , $\gamma$ is said to be 
{\smbf one-way infinite}; 
if~$J$ is bounded and $\gamma$ is simple then $\gamma$ is called an {\smbf arc}. 

The geometric model we use for $\mathbb{H}^2$ is the {\em Poincar\'e disk model}:  
Topologically it is the interior of the compact unit disk $D^2$ in 
the Euclidean plane, the boundary circle $\partial D^2$ is then called 
the {\smbf circle at infinity}. 
The {\smbf complete geodesics} in $\mathbb{H}^2$ are the open diameters, and the open 
arcs of circle orthogonal to $\partial D^2$. The two frontier points of a geodesic 
lie on the circle at infinity and are called its {\smbf endpoints}. 
The complete geodesics in a hyperbolic surface of finite type $\Sigma$ are the 
projections of the complete geodesics in $\mathbb{H}^2$ by the universal covering map. 
A complete geodesic in $\mathbb{H}^2$ 
becomes a curve as soon as a base point, and a tangent vector at it (i.e. a 
{\smbf direction} on the geodesic) have been chosen, the parametrization by 
$\mathbb{R}$ being given by the directed distance from the base point on the 
geodesic.  
A complete geodesic in $\Sigma$ inherits a parametrization from any of its lifts 
to $\mathbb{H}^2$. If the geodesic is compact, any such parametrization 
is a periodic map, hence it can be made into a map from $S^1$ to $\Sigma$, that 
is, the geodesic can be made into a closed curve, and it is then called a 
{\smbf closed geodesic}. A curve, or arc contained in a complete geodesic and 
parametrized by length is called a {\smbf geodesic curve}, or a 
{\smbf geodesic arc}. 

When they belong to a geodesic lamination $\cal L$, the geodesics are called 
the {\smbf leaves} of $\cal L$. 
Likewise, the curves in a topological lamination are called its {\smbf leaves}. 
The complement of $\cal L$ 
is a union of connected components, called the {\smbf complementary regions}
(or {\em principal regions} in~\cite{CB88}) of $\cal L$ in $\Sigma$, and 
it has same area as $\Sigma$, i.e. $\cal L$ has zero area~\cite{CEG87}. 
A complementary region $U$ is the interior of a 
hyperbolic surface with geodesic boundary~\cite[Lem. 4.2.6]{CEG87}.

A leaf of $\cal L$
isolated from the rest of $\cal L$ at least from one side is called
a {\smbf boundary leaf}, and a boundary leaf isolated from the rest of $\cal L$
from both sides is called an {\smbf isolated leaf}. The connected components
of the frontier in $\Sigma$ of the complementary regions all correspond (not
necessarily bijectively) to the boundary leaves~\cite[Lem. 4.2.6]{CEG87}.
A connected component of the frontier of a complementary region contains at least 
one boundary leaf, but it may also contain non-boundary leaves at which boundary 
leaves accumulate (this happens to compact leaves around which isolated leaves 
spiral from both sides), a leaf contained in a connected component of the 
frontier of a complementary region is called a {\smbf frontier leaf}.
The number of complementary regions of $\cal L$ is finite, as is the number
of boundary leaves~\cite[Cor. 4.2.7]{CEG87}.
The set of geodesics obtained by removing from $\cal L$ the (finitely many) 
isolated leaves is still a geodesic lamination ${\cal L}'$, called the 
{\smbf derived} lamination of $\cal L$ in~\cite{CB88}. The derived lamination is 
always a lamination with compact support, which means that no leaf of ${\cal L}'$ 
ends up at a puncture (thus in particular there are only finitely many leaves of 
$\cal L$ ending up at a puncture), moreover ${\cal L}''' = {\cal L}''$ (so that  
the number of frontier leaves of $\cal L$ is finite too), finally, if $\cal L$ 
has no compact leaf then ${\cal L}'' = {\cal L}'$~\cite[Th. 4.2.8]{CEG87}. 

\subsection{Graphs and Train Tracks}
\label{GraphAndTTsec}
\begin{figure}
 \centering
 \includegraphics[width=8cm]{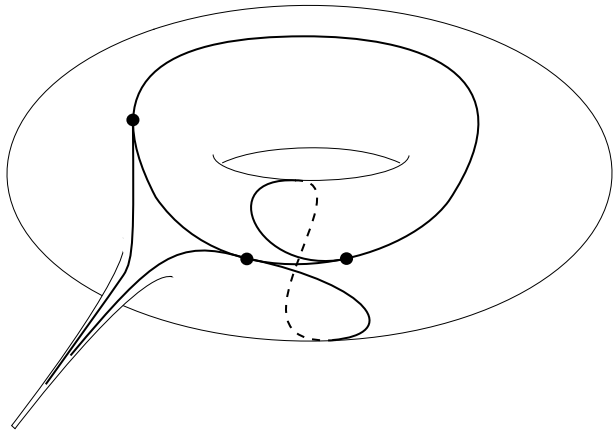}
 \caption{A train track on a once punctured torus with four vertices (one 
is the puncture) and six edges}
 \label{Fig.1}
\end{figure}

A {\smbf graph} is defined by a set of {\smbf vertices}, a set of 
{\smbf edges}, and a map mapping each edge to an unordered pair of (not 
necessarily distinct) vertices; when the pair of vertices of an edge is 
ordered, the edge is said to be {\smbf directed}, the smaller vertex 
is then called its {\smbf origin} and the bigger one its {\smbf end}. 
An {\smbf embedding of a graph} in a surface $\Sigma$ is a pair of 
maps, one sending the vertices into $\Sigma$ union its punctures, and the 
other one sending the edges to simple curves in $\Sigma$, in such a 
way that the endpoints of each curve are the images of the vertices of its 
corresponding edge, and the curves may meet only at their endpoints. 
When a graph is embedded, we identify it with its image by the embedding. 
A {\smbf train track}~\cite{Thu80, PH93} is a finite embedded graph such that: 
The edges are $C^1$ curves, there is a well-defined tangent space at every 
vertex in $\Sigma$ (see Figure~\ref{Fig.1}), and no connected component of 
$\Gamma$'s complement in 
$\Sigma$ is one of the following {\em forbidden regions}: 
\begin{itemize}
				\item a disk whose frontier, as a curve has less than three points 
where it is not 
				$C^1$ (punctures count as such points) 
        \item an annulus whose two frontier components, as curves are closed 
$C^1$ 
        \item a once punctured disk whose frontier component, as a curve is 
closed $C^1$ 
\end{itemize}
\noindent

\noindent
We also make use 
of infinite train tracks in $\mathbb{H}^2$, but they are always preimages by 
the universal covering map of ordinary, that is, finite train tracks on a 
hyperbolic surface of finite type. 

A {\em path} on a graph is a set of directed edges indexed by a subset of 
consecutive integers in $\mathbb{Z}$, in such a way that for every pair 
of edges indexed by a pair $(i, i+1)$, the end of the edge $i$ is the origin 
of the edge $i+1$. 
A path on an embedded graph is also a curve as soon as a base point is chosen, 
the parametrization coming from the 
edge embeddings. On a train track, an {\smbf admissible edge path}, or a 
{\em trainpath} for some authors~\cite{PH93}, is a path everywhere $C^1$ 
as a curve, such that no two consecutive edges share a 
puncture (i.e. admissible paths `do not cross' punctures). From now on, by a 
{\smbf path} on a train track we always mean an admissible edge path. 
A {\smbf connected} train track is a train track of whose any two vertices 
can be joined by a path.  

\subsection{Homotopy}
\label{Homotopy}
An {\smbf ambient isotopy} of a topological space is 
a homotopy $H$ between the identity map and a self-homeo\-morphism $\varphi$ 
of the space such that for each $t \in [0, 1]$, $H(\cdot, t)$ is a self-homeomorphism. 
Given a subset $B$ of the space, $B$ and $\varphi(B)$ are 
then said {\smbf ambient isotopic}. 
On surfaces train tracks and laminations are considered up to ambient isotopy. 
Note that if two embeddings of the same finite graph in a hyperbolic surface of finite 
type $\Sigma$ are isotopic, it is possible to make an ambient isotopy involving 
only a neighborhood of $\Sigma$'s punctures in such a way that both coincide in 
this neighborhood. Then, since the complement of such a neighborhood is a 
compact surface, there is an ambient isotopy deforming one embedding onto 
the other one~\cite[\S 4, Cor.3 and 4]{CVM13}. 
Hence, when we deform a train track onto another one by some isotopy, 
we will just check that the edges are $C^1$ in both train tracks, 
and that the tangent spaces are well-defined and are compatible at the corresponding 
vertices. 

In this paper, two arbitrary curves are said {\smbf homotopic} (resp. 
{\smbf isotopic}) if they admit parametrizations making them homotopic (resp. 
isotopic) using a {\em uniformly continuous} map. This additional condition 
is made necessary by the definitions of {\em carrying}, that we give and 
comment on in the next subsection, and by the fact that the curves we 
consider may not be compact.  
As said before, choosing a base point on a path on a train 
track, resp. a base point and a direction on a geodesic, yields a parametrization 
as a curve for the path, resp. the geodesic. 
Now, a change of base point can be realized by a homotopy, in the above sense, 
between the two parametrized curves, so 
we can speak about homotopy between a curve and a path considered up to index 
translation, and similarly between a curve and a geodesic. 
As for geodesics, note that a curve homotopic to a complete geodesic in a 
surface is necessarily either two-way infinite or closed.

\begin{figure}
 \centering
 \includegraphics[width=8cm]{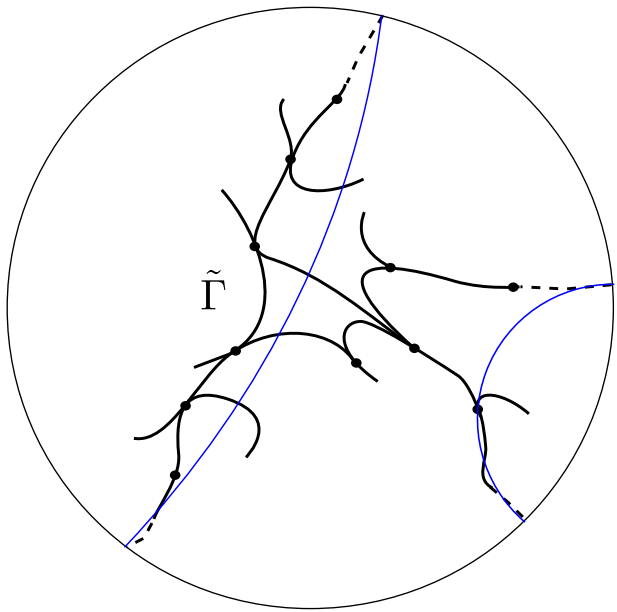}
 \caption{Part of the preimage $\tilde{\Gamma}$ in $\mathbb{H}^2$ (in bold lines) of 
a train track $\Gamma \subset \Sigma$, 
with two complete geodesics (in thin lines) homotopic to paths on $\tilde{\Gamma}$}
 \label{Fig.2}
\end{figure}

Homotopies between curves and paths on a train track 
$\Gamma$ lift to $\mathbb{H}^2$, as homotopies between curves and paths on 
the preimage $\tilde{\Gamma}$ of $\Gamma$ by the universal covering map. 
(See Figure~\ref{Fig.2})

The following lemma is known for compact hyperbolic surfaces with geodesic 
boundary~\cite[Lem. 20]{ZB04}, and can be derived from it.

\begin{Lem}
On a hyperbolic surface of finite type, (1) any path on a train track 
$\Gamma$ lifts to an embedded path in $\tilde{\Gamma}$, (2) any two paths in 
$\tilde{\Gamma}$ which diverge at some vertex never meet again, and (3) every 
path in $\tilde{\Gamma}$ is homotopic to a (unique) geodesic. 
\label{Bona_Lem}
\end{Lem}
\PR
The first two assertions are combinatorial, and will follow easily from 
the construction presented below, only the third one needs a topological 
argument.

Let $\Sigma$ be a non-compact hyperbolic surface of finite type, with a 
train track $\Gamma$ on it. Take a union of pairwise disjoint open 
neighborhoods of $\Sigma$'s punctures in such a way that this union is 
allowed to meet only the non-compact edges of $\Gamma$, and along at most 
one curve for each edge. We further ask that removing these neighborhoods 
results in a (necessarily compact) surface $\Sigma'$ with $C^\infty$ 
boundary, so that the Riemannian metric of $\Sigma$ induces one on 
$\Sigma'$. 
We change the metric on $\Sigma'$ in order to make it into a  
hyperbolic surface with geodesic boundary $\Sigma'_h$. Since 
$\Sigma'$ is compact the two metrics are Lipschitz equivalent. 
Next, in $\Sigma$ we remove from $\Gamma$ the non-compact edges if any, 
and inductively, the vertices where the tangent space has become undefined 
as well as the edges adjacent to them. 
If the resulting graph, necessarily in $\Sigma'_h$, is empty, this implies 
that there are only finitely many two-way infinite paths on $\Gamma$, each 
path being homotopic to a geodesic joining two punctures, and the conclusions 
of the lemma obviously hold. 
If the resulting graph is not empty, it is a disjoint union $\Gamma'$ of  
compact train tracks in $\Sigma'_h$, to each of which~\cite[Lem. 20]{ZB04} 
applies, hence the first two conclusions of the lemma hold for $\Gamma'$ in 
$\Sigma'$. The third one holds by Lipschitz equivalence of the metrics. 
Note that no component of $\partial \Sigma_h$ is carried by $\Gamma'$, 
for if that was the case there would be a forbidden region in the complementary 
set of $\Gamma$, in particular any path carried by the preimage 
$\tilde{\Gamma'}$ of $\Gamma'$ by the universal covering (of $\Sigma$) has 
two distinct endpoints on the circle at infinity. 
Now we add back to $\Sigma'$ the punctured neighborhoods and to $\Gamma'$ 
the removed vertices and edges. Taking into account that these edges involve 
only paths ending up at punctures, and that a path joining a vertex to a puncture 
has finitely many edges, the lemma is proved. 
\EPR

\noindent
It is a consequence of the above lemma, that any two paths on a train track 
which are homotopic differ by an index translation in $\mathbb{Z}$. In other 
words, up to index translation there is only one path homotopic to a given 
curve. 

\subsection{Carrying}
\label{Carrying}
The definitions of carrying do not need the notion of lamination, so 
we state them for sets of curves. 
A set of everywhere smooth curves is {\smbf carried} by a train track $\Gamma$ 
if there is a $C^1$~map from $\Sigma$ to itself homotopic to the identity, 
non-singular on the tangent space of each curve, sending each curve to a path 
on $\Gamma$~\cite[Def. 8.9.1]{Thu80}. 
There are two more notions of carrying, valid for curves not 
necessarily everywhere smooth. 

\begin{figure}
 \centering
 \includegraphics[width=14cm]{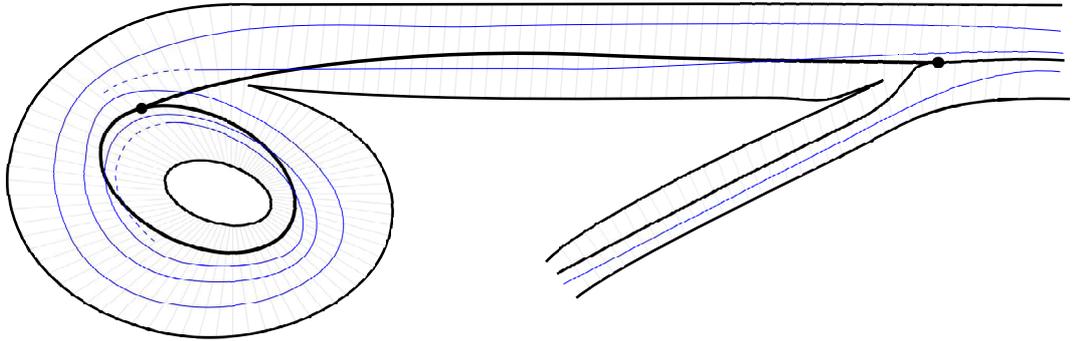}
 \caption{A fibered neighborhood of a train track near two vertices, with 
three pieces of strongly carried leaves}
 \label{Fig.3}
\end{figure}

The first one, already used by Thurston is {\smbf strong carrying} (we use the 
terminology of~\cite{ZB04}). A set of curves $\cal C$ in a surface $\Sigma$ 
is strongly carried by a train track 
$\Gamma$ if there is a regular neighborhood $N \supset {\cal C}$ of $\Gamma$ 
realized as a disjoint union of piecewise $C^1$, or PL arcs, called {\smbf fibers}, 
to which every path on $\Gamma$ is everywhere transverse, and such that each curve 
of ${\cal C}$ 
is also everywhere transverse to the fibers~\cite{PH93, ZB04} (see Figure~\ref{Fig.3}); 
$N$ is called a {\smbf fibered neighborhood} of $\Gamma$. 
Here, since regular neighborhoods and transversality are used, the smooth, or PL 
structures of $\Sigma$ are necessary, and accordingly the carried lamination 
has to be piecewise $C^1$, or PL. 
It is a fundamental result that every geodesic lamination (hence formed of 
everywhere smooth curves) on a hyperbolic surface of finite type has strongly 
carrying train tracks~\cite[Sect. 8.9]{Thu80}, \cite[Th.1.6.5]{PH93}. 
The link between carrying and strong carrying is the following: For every 
lamination carried by a train track $\Gamma$ there is a homeomorphism of $\Sigma$ 
isotopic to the identity map which sends it onto a lamination strongly carried 
by $\Gamma$~\cite{Thu80}, so these two notions of carrying are close.

The second notion of carrying is {\smbf weak carrying}: a set of curves $\cal C$ 
on $\Sigma$ is weakly carried by a train track $\Gamma$ if each curve 
in $\cal C$, is homotopic to a path on $\Gamma$~\cite{ZB04}. 
Weak carrying uses homotopies which are uniformly continuous maps, so there 
is {\it a priori} a dependency on the uniform structure induced by $\Sigma$'s 
metric space structure, however it is known that a change of hyperbolic metric 
on a hyperbolic surface of finite type always induces a Lipschitz equivalent 
metric space structure. 
For weak carrying, as the terminology suggests we have the following: 
\begin{Lem}
Strong carrying implies weak carrying. 
\label{StrongImpliesWeak}
\end{Lem}
\PR
The way curves are deformed onto paths is not a part of the 
definition of strong carrying, however there are natural homotopies: The ones 
obtained by sliding each curve along the transverse fibers towards 
the train track. 
Now, on a hyperbolic surface of finite type there are always regular neighborhoods 
of the train track which can be fibered in such a way that the fibers have uniformly 
bounded lengths, so these natural homotopies can be taken to be uniformly continuous, and 
strong carrying implies weak carrying. 
\EPR

\noindent 
And also: 
\begin{Lem}
Carrying implies weak carrying. 
\label{ImpliesWeak}
\end{Lem}
\PR
A map realizing the carrying by some train track of a lamination made of smooth 
curves can be modified in a neighborhood of the punctures so as to become uniformly 
continuous, so carrying also implies weak carrying. 
\EPR

\noindent
In Definition~\ref{lam_def} by carrying we mean weak carrying, so according to the above 
lemmas, the equivalence between Definitions~\ref{geodlam_def} 
and~\ref{lam_def} holds for the three forms of carrying.

\section{A Topological Lamination is Isotopic to a Geodesic Lamination}
\label{FromTopToGeod}
In this section by `surface' we mean a hyperbolic surface of finite type. 
We begin with a lemma which allows making topological curves into PL 
ones. It proves useful for the subsequent results since it allows the use of  
general position arguments from PL Topology: 

\begin{Lem}
For any simple infinite curve in a surface $\Sigma$, there is an 
isotopy making it a PL curve. 
If the curve is a one-way infinite curve ending up at a puncture 
in a leaf of a topological lamination, the isotopy can be realized 
by an ambient isotopy avoiding the rest of the lamination.
\label{ContinuousToPL_lem}
\end{Lem}

\noindent 
The proof uses methods of classical planar topology, but it is quite long and technical, 
so we postpone it to the Appendix. This lemma has the following corollary:

\begin{Cor}
Let $\gamma$ be a two-way infinite (resp. closed) curve in a 
surface $\Sigma$ homotopic to a two-way infinite (resp. closed) geodesic, 
then $\gamma$ is isotopic to it.
\label{HomotopyAndIsotopy_cor}
\end{Cor}
\PR 
For closed curves, according to~\cite[Th. A1]{Eps66} we can suppose $\gamma$ 
is PL, and since a closed curve homotopic to a closed geodesic is essential,  
\cite[Th. 2.1]{Eps66} gives the desired isotopy. 
Now let $\gamma$ be two-way infinite, and $\gamma_g$ be the geodesic homotopic to it. 
Let $\gamma'$ be a PL 
curve isotopic to $\gamma$ given by Lemma~\ref{ContinuousToPL_lem}, then it 
suffices to prove that $\gamma'$ is isotopic to $\gamma_g$. We lift these 
two curves to $\mathbb{H}^2$ in such a way that the two lifts $\tilde{\gamma'}$ 
and $\tilde{\gamma_g}$ share the same two endpoints on the circle at infinity.
As a first case, if $\gamma'$ and $\gamma_g$ are disjoint, then $\tilde{\gamma'}$ 
and $\tilde{\gamma_g}$ bound an infinite stripe $\mathbb{R} \times [0, 1]$, which 
embeds in $\Sigma$ by the covering projection since 
the image of $\mathbb{R} \times \{0, 1\}$, which is $\gamma' \cup \gamma_g$ does, 
and this provides the desired isotopy. 
The second case is when $\gamma'$ and 
$\gamma_g$ intersect; then we put them in general position, and the intersections 
of $\tilde{\gamma'}$ with the translates of $\tilde{\gamma_g}$ delimit in 
$\mathbb{H}^2$ an enumerable number of disks, any two of which are either 
nested or have disjoint interiors. Then the innermost ones all embed in 
$\Sigma$ by the covering projection, so they can be used to perform inductively 
isotopies in $\Sigma$ until all intersections are erased, making $\gamma'$ and 
$\gamma_g$ disjoint curves, and by the first case above there is an isotopy 
between these two curves.\EPR

\noindent
This corollary suggests that in Definition~\ref{lam_def}, the term ``isotopies'' 
could be replaced by ``homotopies''. But there is a subtlety here: A closed curve 
is parametrized by $S^1$, as well as by $\mathbb{R}$, using a periodic map.  
Then, for instance an embedded two-way infinite curve spiraling along a closed 
geodesic in a compact annulus in $\Sigma$ is homotopic to this geodesic. However 
there can be no isotopy, since one map is an embedding of $\mathbb{R}$ and the 
other one is not. 
This is why in the statement of the following Proposition~\ref{LeafByLeaf_prop}, 
homotopies are used instead of isotopies. Nevertheless, up to a convention that 
we will make after Corollary~\ref{replace_cor}, this proposition   
proves one direction of the equivalence: 

\begin{Prop}
Let $\cal T$ be a topological lamination on a surface $\Sigma$. Then $\cal T$ can 
be deformed onto a geodesic lamination $\cal L$ by curve-by-curve homotopies. 
\label{LeafByLeaf_prop}
\end{Prop}
\PR
Let $\cal T$ be a topological lamination weakly carried by a train track $\Gamma$,  
as in Definition~\ref{lam_def}. We denote by $p: \mathbb{H}^2 \rightarrow \Sigma$ 
the universal covering, and by $\tilde{\Gamma}$ a component of $p^{-1}(\Gamma)$. 
Since each of $\cal T$'s curves has a lift homotopic to one path on 
$\tilde{\Gamma}$, which is homotopic to a geodesic by~Lemma~\ref{Bona_Lem}, 
$\cal T$ can be deformed by curve-by-curve homotopies onto a set ${\cal L}_0$ of 
geodesics. The geodesics in this set are simple and pairwise disjoint, because 
their corresponding curves in $\cal T$ are, and the closure in $\Sigma$ of a set 
of pairwise disjoint simple geodesics is a geodesic 
lamination $\cal L$~\cite[Lem. 3.2]{CB88}. 
Suppose there is a geodesic $\gamma$ in ${\cal L} \setminus {\cal L}_0$, then 
$\gamma$ is the pointwise limit of a sequence $\{\gamma_n\}_{n \in \mathbb{N}}$ 
of geodesics in ${\cal L}_0$. We lift the $\gamma_n$ and $\gamma$ to geodesics 
$\{\tilde{\gamma}_n\}_{n \in \mathbb{N}}$ and $\tilde{\gamma}$ in $\mathbb{H}^2$ 
in such a way that the pairs of endpoints of the $\tilde{\gamma}_n$ lie among 
those of $\tilde{\Gamma}$'s paths and converge 
to the pair of endpoints of $\tilde{\gamma}$, and since each $\tilde{\gamma}_n$ 
has a unique corresponding path on $\tilde{\Gamma}$ with same endpoints, 
$\tilde{\gamma}$ has also an associated path $\tilde{\delta}$ with same 
endpoints. 
Again by Lemma~\ref{Bona_Lem}, $\tilde{\gamma}$ is homotopic to $\tilde{\delta}$, 
hence $\gamma$ is homotopic to the path $\delta = p(\tilde{\delta})$. 

If each curve in $\cal T$, closed or two-way infinite, has its corresponding geodesic 
in $\cal L$ accordingly closed or two-way infinite, we apply 
Corollary~\ref{HomotopyAndIsotopy_cor} to deduce that ${\cal L}_0$ is isotopic to 
$\cal T$ by leaf-by-leaf isotopies. Now, $\Gamma$ carries ${\cal L}_0 \cup \gamma$, but 
this contradicts the fact that $\cal T$ should remain maximal with 
respect to being carried by $\Gamma$, whatever the isotopies applied to its curves, 
hence there is no such geodesic as $\gamma$, so ${\cal L}_0 = \cal L$. 

If there are two-way infinite curves in $\cal T$ homotopic to closed geodesics 
(necessarily finite in number) we denote these geodesics by $\gamma_i, i=1, \cdots, k$, 
and we consider $k$ pairwise disjoint closed annulus neighborhoods of the $\gamma_i$. 
We first make an isotopy in each annulus in order to move all the curves of 
${\cal L}_0 \cup \gamma$ other than the $\gamma_i$ away from them, 
without creating intersections among the curves. Next, in the annuli we replace 
the $\gamma_i$ by homotopic spiraling two-way infinite simple curves; we obtain in 
this way a set of curves ${\cal L}_1 \cup \gamma'$. Now, 
according to Corollary~\ref{HomotopyAndIsotopy_cor}, 
these spiraling curves are isotopic to their corresponding curves in $\cal T$, 
so that ${\cal L}_1$ is isotopic curve by curve to $\cal T$. 
However, the set of pairwise disjoint pairwise non-homotopic curves 
${\cal L}_1 \cup \gamma'$ is carried by $\Gamma$, again 
contradicting maximality of $\cal T$'s carrying. 
Hence here too ${\cal L}_0 = {\cal L}$. 
\EPR

\noindent
The key-fact with this proposition, is that it provides a natural bijection between 
the leaves of a given topological lamination $\cal T$ in $\Sigma$ and those of a 
unique geodesic lamination 
$\cal L$, and this allows defining the {\smbf boundary leaves}, the {\smbf isolated 
leaves} and the {\smbf frontier leaves} of $\cal T$ as the leaves homotopic to the 
boundary leaves, isolated leaves, and frontier leaves respectively, of $\cal L$. 
The frontier leaves of $\cal T$ lift through $p$, the universal covering map, to 
$\mathbb{H}^2$, and (using the Schoenflies Theorem) they delimit disks in $\mathbb{H}^2$ 
union the circle at infinity such that their interiors contain no lift of $\cal T$'s 
leaves. This leads to:

\begin{Def}
Let $\cal T$ be a topological lamination on a surface $\Sigma$, and $\tilde{E}$ be an 
open disk in $\mathbb{H}^2$ bounded by lifts of $\cal T$'s frontier leaves such that 
$\tilde{E}$ contains no lift of a leaf of $\cal T$. 
The image of $\tilde{E}$ by $p$ is called a {\smbf complementary region} of $\cal T$. 
\label{CompReg_def}
\end{Def}

\noindent
The complementary regions are open sets, finite in number, but unlike for geodesic 
laminations their union need not 
be equal to $\Sigma \setminus \cal T$ (it may even not be a full measure set 
in $\Sigma$). 

\begin{figure}
 \centering
 \includegraphics[width=8cm]{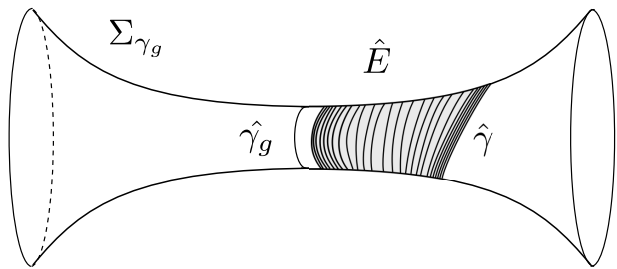}
 \caption{The covering $\Sigma_{\gamma_g}$}
 \label{Fig.4}
\end{figure}

\begin{Lem}
Let $\cal T$ be a topological lamination on a surface $\Sigma$. If a two-way infinite 
leaf $\gamma$ of $\cal T$ 
is homotopic to a simple closed geodesic $\gamma_g$, there is a complementary region $E$ 
of $\cal T$ such that $E \cup \gamma$ contains a compact annulus whose boundary 
components are PL curves homotopic to $\gamma$. 
\label{avoid1_lem}
\end{Lem}
\PR
Note that a closed leaf in a geodesic lamination is always a frontier leaf, so $\gamma_g$ 
is in the geodesic lamination corresponding to $\cal T$ given by 
Proposition~\ref{LeafByLeaf_prop}, thus $\gamma$ is a frontier leaf for $\cal T$. 
The fundamental group of $\Sigma$ is considered at a base point $x \in \gamma_g$. 
Topologically, the cyclic covering $\Sigma_{\gamma_g}$ corresponding to the homotopy 
class of $\gamma_g$ is an open annulus, whose fundamental group $G_c \cong \mathbb{Z}$ 
is generated by the homotopy class of a lift $\hat{\gamma_g}$ of $\gamma_g$ at some 
point $\hat{x}$ over $x$. From now on, we denote the elements 
of $G_c$ by the integers in $\mathbb{Z}$ they correspond to. 
Lifting the homotopy between $\gamma$ and $\gamma_g$ to $\Sigma_{\gamma_g}$ gives a 
curve $\hat{\gamma}$ in $\Sigma_{\gamma_g}$, that we lift further to $\mathbb{H}^2$ to 
obtain a family of curves $\{\tilde{\gamma_i}\}_{i \in G_c}$, all sharing the same two 
endpoints on the circle at infinity, and all pairwise disjoint, because $\hat{\gamma}$ 
is simple. 
Each pair of consecutive $\tilde{\gamma_i}$, 
$(\tilde{\gamma_i}, \widetilde{\gamma_{i+1}})$ 
bounds an open disk $\tilde{E_i}$, and $\tilde{E_i}$ 
projects via the universal covering to a connected open set $\hat{E}$ in 
$\Sigma_{\gamma_g}$ (shaded in Figure~\ref{Fig.4}), and to another one $E$ in 
$\Sigma$, so $E$ is a complementary region. 

\noindent
For $k > 0$ we denote by 
$\tilde{F_k}$ the set $\bigcup_{-k \le i \le k} (\tilde{\gamma_i} \cup \tilde{E_i}) \setminus 
\widetilde{\gamma_{-k}}$. The $\tilde{F_k}$ form an exhaustive sequence 
of open disks in $\tilde{F} = \bigcup_{i \in G_c} (\tilde{\gamma_i} \cup \tilde{E_i})$, 
thus in particular $\tilde{F}$ is an open set. 
As a consequence, $\hat{F} = \hat{E} \cup \hat{\gamma}$ is an open set in $\Sigma_{\gamma_g}$, and 
so is $F = E \cup \gamma$ in $\Sigma$. 
We claim that the inclusion of $\hat{F}$ in $\Sigma_{\gamma_g}$ induces an 
injective morphism on fundamental groups; to see it take a simple closed curve 
$\hat{\delta}$ in $\hat{F}$ bounding a disk in $\Sigma_{\gamma_g}$; $\hat{\delta}$ 
lifts to a closed curve $\tilde{\delta}$ in $\mathbb{H}^2$ also bounding a disk, 
but that curve lies in some $\tilde{F_k}$ because $\tilde{\delta}$ is compact. Now, 
$\tilde{F_k}$ is a disk, hence the disk bounded by $\tilde{\delta}$ lies in 
$\tilde{F_k} \subset \tilde{F}$ and as a consequence, the disk bounded by 
$\hat{\delta}$ lies in $\hat{F}$, proving the claim. This claim implies that $\hat{F}$ 
is either a disk or an annulus. 
We prove now that it is not a disk: We take a point $\tilde{y}$ in $\tilde{E_i}$, 
with image $\hat{y}$ in $\Sigma_{\gamma_g}$. The lift to $\mathbb{H}^2$ at $\tilde{y}$ 
of a loop based at $\hat{y}$ and freely homotopic to $\gamma_g$ gives a path whose 
other endpoint $\tilde{y'}$ lies in $\widetilde{E_{i+1}}$. We link $\tilde{y}$ and 
$\tilde{y'}$ by a path $\tilde{\eta}$ in the disk 
$\tilde{E_i} \cup \widetilde{\gamma_{i+1}} \cup \widetilde{E_{i+1}}$. Note that 
$\tilde{\eta}$ projects to a loop $\hat{\eta}$ in $\hat{F}$ which is not null-homotopic. 
Hence $\hat{F}$ is an open annulus, 
and as such, it contains a compact annulus with PL boundary, which projects to $\Sigma$, 
as a subset of $F = E \cup \gamma$. 
Since the projection of the boundary of the above annulus can be seen as a pair of PL closed 
curves homotopic to $\gamma_g$, which is simple, we can erase their self-intersections if any, 
and also the intersections they may have with each other; moreover these erasures can be 
performed in $F$, so they cobound in $F$ an annulus whose boundary components (as curves) are 
homotopic to $\gamma_g$. 
\EPR

\noindent
The following corollary follows: 

\begin{Cor}
Let $\cal T$ be a topological lamination on a surface $\Sigma$ containing a non-compact 
leaf $\gamma$ homotopic to a simple closed curve. Then there is a lamination 
${\cal T}_\delta$ sharing with $\cal T$ the same leaves except $\gamma$, which has 
been replaced by a simple closed curve $\delta$ homotopic to $\gamma$.
\label{replace_cor}
\end{Cor}
\PR
The fact that $\cal T$ is carried by a train track implies that the closed curve $\gamma$ is 
homotopic to is essential in $\Sigma$, hence is freely homotopic to a simple closed geodesic. 
Under these conditions, Lemma~\ref{avoid1_lem} says that there is in the complement of 
${\cal T} \setminus \gamma$ in $\Sigma$, an annulus containing a simple closed curve 
$\delta$ homotopic to $\gamma$. 
Hence ${\cal T_\delta} = ({\cal T} \setminus \gamma) \cup \delta$ is a lamination as desired . 
\EPR

\noindent
Corollary~\ref{replace_cor} allows the following convention: 
\begin{Convention}
In a topological lamination we 
always replace a two-way infinite curve homotopic to a closed one avoiding the rest of the 
lamination, by the closed one. 
\label{ClosedLeaf_conv}
\end{Convention}

\noindent
Note that for geodesics this point
of view is natural, since on a hyperbolic surface, there can be only one
geodesic homotopic to a curve and it is closed if the curve is closed. With this 
convention, the combination of Corollary~\ref{HomotopyAndIsotopy_cor} and 
Proposition~\ref{LeafByLeaf_prop} gives one direction of the equivalence between 
Definitions~\ref{geodlam_def} and~\ref{lam_def}. 

\section{A Geodesic Lamination is a Topological Lamination}
\label{DecGeoLam}

It suffices to prove it for connected geodesic laminations, and in this section 
$\cal L$ denotes such a lamination. We begin by remarking 
that $\cal L$  contains no leaf joining two punctures. This can be seen as follows: 
Suppose $\cal L$ contains such a (necessarily isolated) leaf $\ell$. 
We can write $\cal L = ({\cal L} \setminus \ell) \cup \ell$; however a geodesic 
joining two punctures is a closed set, and so is ${\cal L} \setminus \ell$, so $\cal L$ 
is not connected. 

\begin{figure}
 \centering
 \includegraphics[width=4cm]{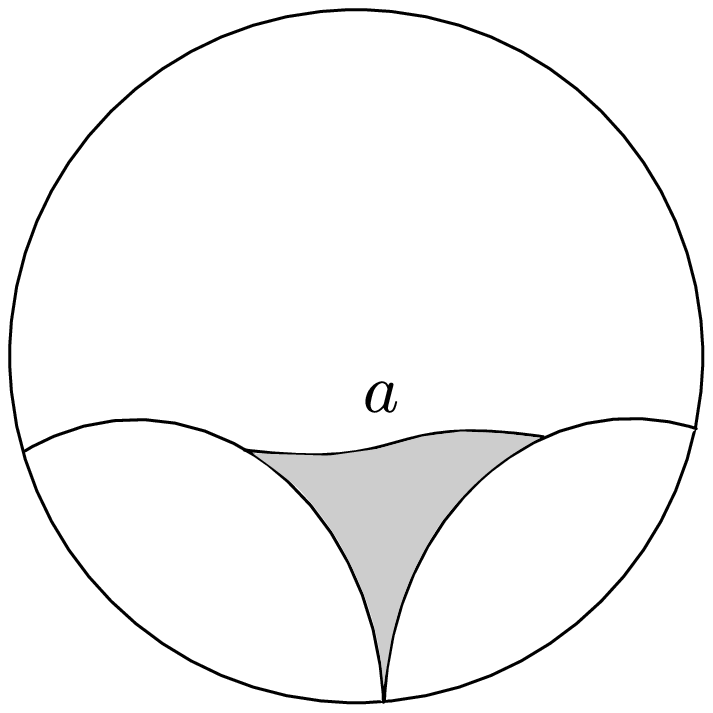}
 \caption{A spike in $\mathbb{H}^2$, with interior shaded}
 \label{Fig.5}
\end{figure}

\begin{figure}
\centering
\includegraphics[width=10cm]{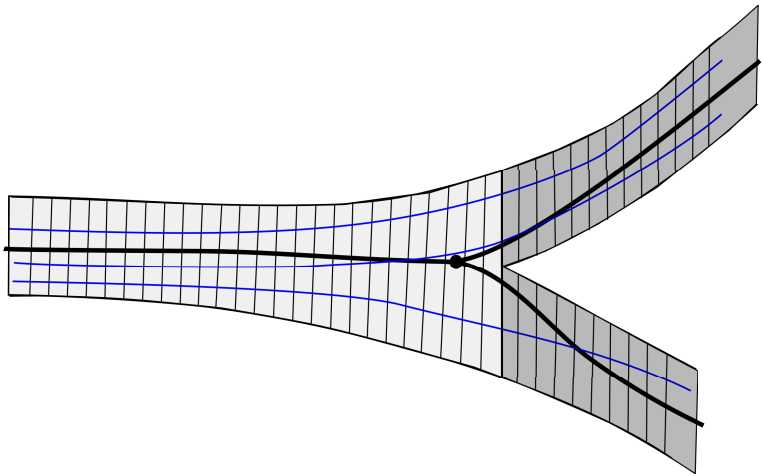}
\caption{A part of $N_1$, formed by the union of three boxes.} 
\label{Fig.6}
\end{figure}

To describe how geodesics in a complementary region of a geodesic lamination 
behave we need the notion of {\em spike}. A {\smbf spike}  is a surface isometric 
to the connected closed region of bounded area in $\mathbb{H}^2$, delimited by two 
disjoint complete geodesics sharing an endpoint, and an arc $a$ meeting each of 
them along one of its endpoints. The arc $a$ is called the {\smbf base} of the spike. 
An {\smbf open spike} is the interior of a spike union the interior of its base. 
(See Figure~\ref{Fig.5})
To work in fibered neighborhoods of train tracks we also need the notion of 
{\em box}. 
A {\smbf box} $D$ in a fibered neighborhood of a train track 
is a union of subarcs of the fibers which forms an embedded closed disk, such 
that its frontier $\partial D$ consists of a pair of disjoint subarcs of two 
(not necessarily distinct) fibers, the {\smbf small sides} of $D$,  
union a pair of piecewise $C^1$ disjoint arcs everywhere transverse to the 
fibers, the {\smbf large sides} of $D$ (see Figure~\ref{Fig.6}). 

Let $N_\varepsilon({\cal L})$ be the set of points of 
$\Sigma$ at distance no more than $\varepsilon$ from $\cal L$. Then
$N_\varepsilon({\cal L})$ is the union (in a non-unique way) of a 
compact topological surface embedded in $\Sigma$ and finitely many 
neighborhoods of spikes. It is known that $N_\varepsilon({\cal L})$ can be 
decomposed as a disjoint union of piecewise $C^1$ arcs making 
$N_\varepsilon({\cal L})$ a fibered neighborhood of 
$\cal L$~\cite{Thu80, PH93}.  
Also, collapsing each fiber to one point yields a corresponding graph 
$\Gamma(\varepsilon, {\cal L})$, 
that embeds in $\Sigma$ with edges $C^1$-embedded in such a way that the 
tangent space at each vertex is well-defined; moreover if $\varepsilon$ is 
small enough the connected components of 
$\Sigma \setminus \Gamma(\varepsilon, {\cal L})$ satisfy the conditions in 
order that $\Gamma(\varepsilon, {\cal L})$ is a train 
track~\cite[Prop. 8.9.2]{Thu80}. 
The boundary components of $N_\varepsilon({\cal L})$ are 
differentiable everywhere except at finitely many points called {\smbf cusps}.
Adequately stretching the fibers in small neighborhoods of cusps makes 
$N_\varepsilon(\cal L)$ into a new fibered neighborhood $N_1$ of 
$\cal L$ such that no two cusps lie in the same fiber; we call the 
corresponding graph obtained by collapsing each fiber to one 
point, $\Gamma_1$. After ambient isotopy the graph $\Gamma_1$ can be embedded 
in such a way that $N_1$ is a regular neighborhood of it, and that $\Gamma_1$ 
strongly carries $\cal L$. Note that as a graph, $\Gamma_1$ is trivalent. 

Suppose there exists a simple geodesic $\gamma_g$ in the complement of 
$\cal L$ homotopic to a path on $\Gamma_1$. We have: 

\begin{Lem}
The geodesic $\gamma_g$ decomposes as a geodesic arc union two disjoint 
geodesic one-way infinite curves $\gamma_{g}^-$ and $\gamma_{g}^+$ in $N_1$, 
each contained in a spike bounded by boundary leaves of $\cal L$. 
Accordingly, the path on $\Gamma_1$ homotopic to $\gamma_g$ decomposes as 
a finite path union two one-way infinite subpaths of paths homotopic to 
boundary leaves of $\cal L$. 
\label{geod_in_spike_lem}
\end{Lem}
\PR
The first assertion is a straightforward consequence of the structure theorem 
about geodesic laminations~\cite[Th. 4.2.8]{CEG87}: Whatever the direction 
chosen on $\gamma_g$, $\gamma_g$ eventually enters an open spike bounded by 
two boundary leaves of $\cal L$. Moreover, when these spikes become 
sufficiently thin (i.e. the pieces of fibers of $N_1$ joining the two 
boundary leaves become sufficiently short), they are subsets of 
$N_\varepsilon({\cal L}) \subset N_1$, because of $N_\varepsilon({\cal L})$'s 
definition. Hence $\gamma_{g}^-$ and $\gamma_{g}^+$ can be taken  
to be subsets of $N_1$. 
The second assertion is an obvious consequence of the first one. 
\EPR

There is no reason $\gamma_g$ should be entirely contained in 
$N_1$. However: 

\begin{Lem}
The path $\gamma$ on $\Gamma_1$ homotopic to $\gamma_g$ can be moved by a 
homotopy in such a way that it becomes a simple curve $\gamma \subset N_1$, 
disjoint from $\cal L$, and transverse to the fibers of $N_1$. 
\label{Carrying_Lem}
\end{Lem}
\PR
The neighborhood $N_1$ can be seen as a union of boxes, one 
box for each edge of $\Gamma_1$. 

We fix in $\gamma_g$ two one-way infinite geodesic curves $\gamma_{g}^-$ and 
$\gamma_{g}^+$, whose existence is asserted in Lemma~\ref{geod_in_spike_lem}. 
Up to taking them shorter we can assume that $\gamma_{g}^-$ and $\gamma_{g}^+$ 
have their endpoints on sides of two boxes $R^-$ and $R^+$ 
respectively, in such a way that they can be prolonged through $R^-$ and $R^+$. 
We now prolong, say $\gamma_{g}^-$, through $R^-$ using a geodesic arc disjoint 
from the boundary leaves of $\cal L$. 
We carry on this process of prolonging the path without creating intersections 
with $\cal L$'s boundary leaves or self-intersections. If at some of the 
boxes forming $N_1$ we had to intersect $\cal L$'s boundary leaves 
or the already drawn arc to be able to continue, lifting to 
$\mathbb{H}^2$ and applying Lemma~\ref{Bona_Lem} would give a non-erasable 
intersection between a lift of $\gamma$ and a lift of a path homotopic to some 
leaf of $\cal L$, or a non-erasable intersection between two lifts of $\gamma$. 
This would contradict the hypothesis that $\gamma_g$ is simple and misses $\cal L$. 
Hence we can continue, following the finite path whose existence is 
asserted in Lemma~\ref{geod_in_spike_lem}. This being done, we meet $R^+$ in such a way 
that we can join the already drawn arc to $\gamma_{g}^+$, and this can be done 
without creating intersections, again thanks to Lemma~\ref{Bona_Lem}. 
\EPR

We can prove now:

\begin{Prop}
For any connected geodesic lamination $\cal L$, there is a train track carrying it 
in a maximal way. 
\label{Carrying_Prop}
\end{Prop}
\PR
We begin with the neighborhood $N_1$, which, as said in the preceding lemma, 
can be seen as a union of boxes, one box for each edge of 
$\Gamma_1$. 
Suppose $\gamma_g$ is a geodesic in $\Sigma \setminus {\cal L}$ homotopic to 
a path on $\Gamma_1$. We assume first that $\gamma_g$ does not end up at a 
puncture. We denote by $\gamma$ a simple curve in $N_1$ 
homotopic to $\gamma_g$ whose existence is proved in Lemma~\ref{Carrying_Lem}. 

We denote by $\gamma_{g, 1}^-$ and $\gamma_{g, 2}^-$ two geodesics in $\cal L$ 
forming a spike meeting $\gamma$ along a one-way infinite curve. Since 
the leaves of $\cal L$ and $\gamma$ are transverse to the fibers of 
$N_1$, when we go sufficiently far into the spike we can find 
in some fiber a subarc $a$ joining $\gamma_{g, 1}^-$ and $\gamma_{g, 2}^-$, 
and meeting $\gamma$ at one point. 
\begin{figure}
\centering
\includegraphics[width=14cm]{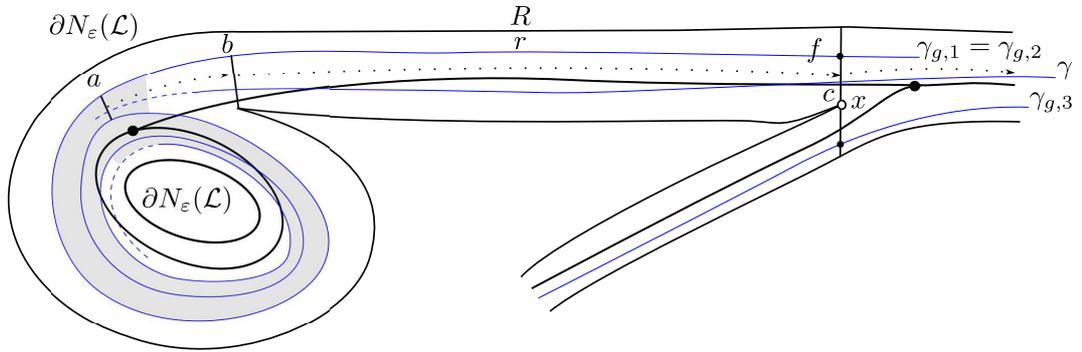}
\caption{A part of $N_1$, where an open spike containing $a$'s interior is shaded. 
In this example $\gamma_{g, 1}^- = \gamma_{g, 2}^-$ spirals around 
a closed geodesic (not represented).} 
\label{Fig.7}
\end{figure} 

As a first step, we translate $a$ along $\gamma$ so as to move away from the 
spike. The translation is possible until the moment we hit a cusp, which 
is also the intersection point between two small sides of two boxes 
of $N_1$. 
This eventually happens since none of $\gamma_{g, 1}^-$ or $\gamma_{g, 2}^-$ 
is homotopic to $\gamma$, and this ends step~1. 
Two boxes share the hit cusp, and we choose the one, $R$, visited by the 
prolonging of $\gamma$. Then, one of $\gamma_{g, 1}^-$ or $\gamma_{g, 2}^-$, 
say the latter, enters also $R$ together with $\gamma$. 
In $R$ there is a thinner box $r$ whose one large side is a 
side of $R$ and the other one is a subarc of $\gamma_{g, 2}^-$, its small sides 
being subarcs of $R$'s small sides. 

\begin{figure}
\centering
\includegraphics[width=10cm]{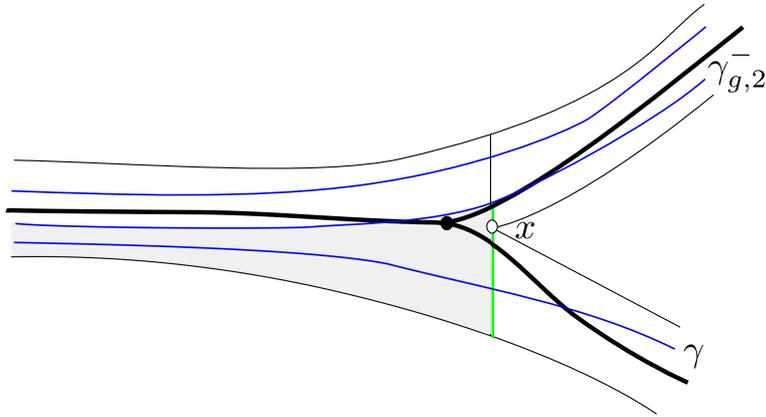}
\caption{Stopping the translation at a hit with a cusp $x$ in the interior of a 
box small side.} 
\label{Fig.8}
\end{figure}

As a second 
step we translate further through $r$ along $\gamma$, and even further 
than $r$ provided we do not have to stop at a cusp $x$. 
However we have to stop some time, ending thus step~2, since, again, $\gamma$ is 
not homotopic to $\gamma_{g, 2}^-$. 
Then, the hit cusp $x$ must be an endpoint of a small side of the 
box resulting from the translation, in some fiber $f$ of $N_1$, as 
shown in Figure~\ref{Fig.7}. 
Indeed, if the hit cusp was in the interior of the small side where we 
stopped translating (see Figure~\ref{Fig.8}), one of the boxes 
forming $N_1$ would not be visited by any leaf of $\cal L$, which is 
impossible. 

\noindent
When this happens, the fiber $f$ necessarily intersects pieces of leaves of 
$\cal L$, and by compactness of $\cal L$ there is a subarc $c \subset f$ 
joining $\gamma_{g, 2}^-$ to another boundary leaf $\gamma_{g, 3}^-$, such 
that $c$ contains no hit with any leaf of $\cal L$ other than 
$\gamma_{g, 2}^-$ and $\gamma_{g, 3}^-$. 

We begin then Step~3; we are in a situation similar to the one we had at 
the beginning of step~1, when we started moving $a$ away from the spike, 
with now the subarc $c$,  that we can translate following $\gamma$ until 
we hit again a cusp, putting an end to this step. 
If the cusp $y$ we eventually hit lies between the hit of 
$\gamma_{g, 2}^-$ and the one of $\gamma$ in the fiber containing $y$,  
there exists in $N_1$ an arc $\delta$ transverse to the fibers, 
joining $x$ and $y$, and disjoint from $\cal L$; $\delta$ always exists 
because no two cusps lie on the same fiber. If $y$ does not lie between 
the hit of $\gamma_{g, 2}^-$ and the one of $\gamma$, we are in the 
same situation as at the beginning of step~2, and we can proceed step~4 in 
a way similar to step~2. 

If this process never stops, then $\gamma_{g, 2}^-$ is homotopic 
to $\gamma$, and this contradicts the hypothesis. 
Hence for some $k \ge 1$, at step $2k+1$ we reach a  
situation where the hit cusp $y$ lies between the hit of $\gamma_{g, 2}^-$ 
and that of $\gamma$ in the fiber containing $y$. We conclude that there 
is always an arc $\delta$ transverse to the fibers, missing $\cal L$, and 
joining $x$ to another cusp $y$. 

We slit open $N_1$ along $\delta$, i.e. we remove from $N_1$ the interior 
of a thin box neighborhood of $\delta$ missing $\cal L$, 
to get a fibered neighborhood $N_2$ containing $\cal L$, in such a way that 
$\cal L$ is transverse to the fibers. After smoothing, $\partial N_2$ has 
two fewer cusps than $N_1$. 

We repeat on $N_2$ the process undergone for $N_1$, and so on. Since the number 
of cusps decreases, we eventually reach a situation 
where, after $n$ slitting open operations, the interior of the obtained fibered 
neighborhood $N_n$ of $\cal L$ contains no more arcs transverse to all the fibers,  
joining two punctures, and missing $\cal L$. 
If the corresponding train track $\Gamma_n$ had a path homotopic to some 
geodesic in $\Sigma \setminus {\cal L}$, thanks to Lemmas~\ref{geod_in_spike_lem} 
and~\ref{Carrying_Lem} we would obtain a simple curve $\gamma$ in $N_n$ whose 
existence would imply the existence of an arc as above, and this is impossible. 
Hence the carrying of $\cal L$ by $\Gamma_n$ is maximal. 

When $\gamma_g$ ends up at a puncture (necessarily only) in one direction, the 
situation is easier: In the direction where the associated path $\gamma$ 
given by Lemma~\ref{Carrying_Lem} does not end up at a puncture it enters a spike 
bounded by two boundary leaves $\gamma_{g, 1}^-$ and $\gamma_{g, 2}^-$ of 
$\cal L$, so the beginning of the above argument applies, and in the other 
direction $\gamma$ runs only along finitely many edges. 
\EPR

\noindent
Proposition~\ref{Carrying_Prop} ends the proof of the equivalence of 
Definitions~\ref{geodlam_def} and~\ref{lam_def} up to isotopy. 

\section{Moving Topological Laminations into their Carrier Train Track Neighborhoods}
\label{Engulf}

In this section we relax the dependency on the embedding surface, by proving 
that a weakly carried lamination satisfies, up to ambient isotopy, 
the first requirement for being strongly carried, that is, to be contained in a regular 
neighborhood of the carrying train track. 

We begin with a lemma about the leaves ending up at a puncture: 

\begin{Lem}
Let $\cal T$ be a topological lamination containing curves ending up at 
punctures, then there is an ambient isotopy sending $\cal T$ to a topological 
lamination such that all the leaves ending up at punctures are PL in a neighborhood 
of the punctures. 
\label{PunctPL_lem}
\end{Lem}
\PR
The set of leaves ending up at a puncture is finite. For each puncture and each leaf 
ending up at it, we choose a one-way infinite curve ending up at that puncture in that 
leaf.
According to Lemma~\ref{ContinuousToPL_lem} there is an ambient isotopy bringing 
each of these one-way infinite curves to a PL curve while fixing the rest of the 
lamination. Composing these finitely many ambient isotopies makes all the leaves 
ending up at punctures, PL near the punctures. 
\EPR 

\begin{figure}
 \centering
 \includegraphics[width=8cm]{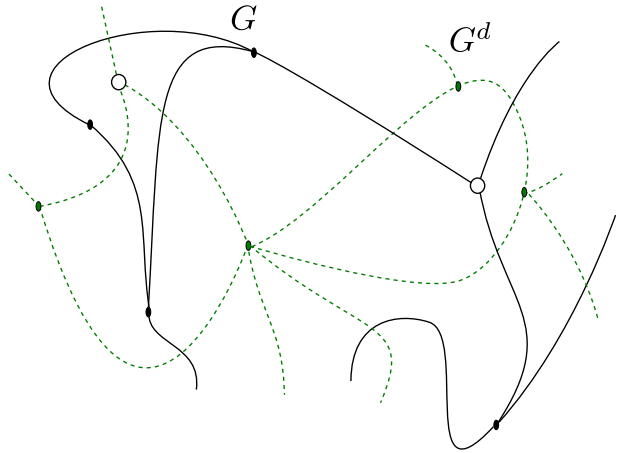}
 \caption{Part of a graph $G$ and of its dual $G^d$. Solid lines represent $G$'s 
edges, narrow dashed ones, $G^d$'s edges, and small circles, punctures of $\Sigma$ 
(there are two of them).}
 \label{Fig.9}
\end{figure}

An embedded graph in $\Sigma$ (not necessarily a train track) {\smbf fills} $\Sigma$ 
if each of the connected components of the complement of the graph is either a disk 
or a once punctured disk. For a filling graph $G$, its {\smbf dual graph} in $\Sigma$ 
is defined as follows: 
In each disk component of $\Sigma \setminus G$ we choose a point. 
These points and the punctures in the once punctured disk components of 
$\Sigma \setminus G$ are joined with each other by PL simple curves, actually 
arcs or one-way infinite curves with pairwise 
disjoint interiors in $\Sigma$, in such a way that each pair of (not necessarily 
distinct) components sharing an edge of $\Gamma$ are joined, and the joining 
curve meets $\Gamma$ in only one point, in that edge. 
The resulting embedded graph $G^d$ is called a graph {\smbf dual to} $G$. 
(See Figure~\ref{Fig.9})

\noindent
Since all the components of $\Sigma \setminus G$ are disks or once punctured disks, 
any two dual graphs are isotopic, hence ambient isotopic, so $G^d$ is essentially 
unique. We have:
\begin{Rem}
Let $G$ be a filling graph in $\Sigma$. Then $G$ is the dual graph of $G^d$. 
\label{Fill_dual_reflexive_rem}
\end{Rem}
\PR
By construction, the connected components of $\Sigma \setminus G^d$ satisfy the same 
topological constraints as those of $\Sigma \setminus G$, with the difference that 
a disk component of $\Sigma \setminus G^d$ contains a preferred point in it, 
that is, the (unique) vertex of $G$ whose adjacent edges hit the frontier of this 
component. And the situation is similar for the once punctured disks of 
$\Sigma \setminus G^d$. Thus we can take these points and punctures as the vertices 
of the dual of $G^d$, these are actually exactly the vertices of $G$. For 
the edges of $G^d$'s dual we just take the edges of $G$. 
\EPR

The following corollary is the core of the proof of Lemma~\ref{Engulfing_fill_lem}: 

\begin{Cor}
Let $G$ be a filling graph in $\Sigma$. Then there is a set of neighborhoods of 
$G^d$'s vertices in $\Sigma$ such that the complement of the 
union of these neighborhood interiors is a regular neighborhood of $G$. 
\label{Engulfing_fill_cor}
\end{Cor}
\PR
This type of result is classical in PL Topology. 
We take disjoint closed PL neighborhoods missing $G$ of each of $G^d$'s vertices, 
each being a PL disk or once punctured disk according to whether the considered vertex 
lies in $\Sigma$ or is a puncture, in such a way that if the frontier of the 
component of $\Sigma \setminus G$ the vertex lies in contains punctures, 
the frontier of the neighborhood contains the same punctures. 
The complement in $\Sigma$ of the union of these neighborhood interiors is a 
PL subsurface of $\Sigma$, which collapses onto $G$, so it is a regular 
neighborhood of $G$. 
\EPR

We need one more graph associated to $G$, which is built using both $G$ and $G^d$. 
A {\smbf puncture graph of} an embedded graph $G$ is a graph $\Lambda_G$, 
whose vertices are the punctures and those vertices of $G^d$ whose neighborhood 
frontiers contain punctures, and the edges are arcs whose interiors lie 
in $\Sigma \setminus G$, joining each such vertex to the punctures on its neighborhood 
frontier. The edges of $\Lambda_\Gamma$ can be taken PL since by 
Lemma~\ref{PunctPL_lem} the frontier leaves are PL near the punctures. 
Again, up to ambient isotopy $\Lambda_G$ is unique. 

\begin{figure}
 \centering
 \includegraphics[width=8cm]{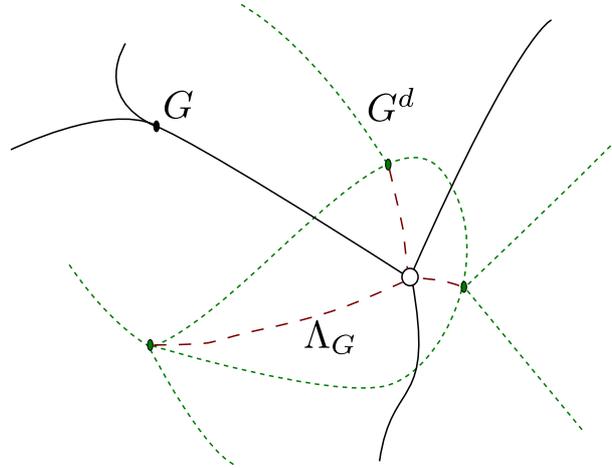}
 \caption{part of a puncture graph $\Lambda_G$. Its edges are represented as wide dashed 
lines}
 \label{Fig.10}
\end{figure} 

\noindent
For connected train tracks we have:

\begin{Lem}
Let $\Gamma$ be a connected train track. 
The vertices of $\Gamma$ which are not punctures all lie in the same connected component 
of the complement of $\Gamma$'s puncture graph $\Lambda_\Gamma$. 
\label{puncture-graph_lem}
\end{Lem}
\PR
Let $\Lambda_\Gamma$ be the puncture graph of $\Gamma$ in $\Sigma$. 
Any path joining two vertices of $\Gamma$ lying in two distinct components of 
$\Sigma \setminus \Lambda_\Gamma$ would have to contain a puncture vertex in its interior, 
so it is would not be admissible, but this contradicts the hypothesis that $\Gamma$ 
is connected. 
\EPR

Let $\cal T$ be a topological lamination carried by a train track $\Gamma$. 
To each frontier leaf of $\cal T$ there corresponds one two-way infinite path on $\Gamma$, 
so to each complementary region $E$ of $\cal T$ we associate the set of 
paths homotopic to its frontier leaves, and call these paths the {\smbf frontier paths  
for} $E$. We lift the frontier paths for $E$ to $\mathbb{H}^2$ through $p$, in such a way 
that their endpoints coincide with the endpoints of the frontier components of a disk 
$\tilde{E} \subset \mathbb{H}^2$ over $E$. Then, the union of these frontier paths' 
lifts delimit regions, lying over connected components of $\Sigma \setminus \Gamma$ that  
we call {\smbf components of} $\Sigma \setminus \Gamma$ {\smbf associated to} $E$, since the 
set they form does not depend on the choice of $\tilde{E}$ over $E$. 
We have:

\begin{Lem}
Let $\cal T$ be a topological lamination in $\Sigma$ weakly carried by a connected  
train track $\Gamma$ which fills $\Sigma$. Then $\Gamma^d$ is isotopic to an embedded 
graph $\Delta$ having one vertex in each component of $\Sigma \setminus \Gamma$ associated 
to $E$, where $E$ runs through the complementary regions of $\cal T$, in such a way that 
$\Lambda_\Delta \subset (\Sigma \setminus {\cal T})$. 
\label{Engulfing_prepa_lem}
\end{Lem}
\PR
Given a complementary region $E$, we choose in it one point for each disk 
component $e$ of $\Sigma \setminus \Gamma$ associated to $E$, if $e \cap E \neq \emptyset$ 
the point for $e$ is taken in $e \cap E$. 
As for the once punctured disk components, we take their punctures. 
We mark each point according to the vertex of the dual graph $\Gamma^d$ it corresponds 
to, including the puncture vertices. We denote by $V_E$ the set thus obtained, 
and by $V$ the union of the $V_E$ where $E$ runs through the complementary regions of 
$\cal T$. 
For each $E$ we slide each non-puncture vertex of the puncture graph of $\Gamma$, 
$\Lambda_\Gamma$, onto the point of 
$V_E$ it corresponds to, in such a way that the edges of $\Lambda_\Gamma$ can be moved 
rel endpoints so as to lie in the union of ${\cal T}'s$ complementary regions.
These isotopies can be obtained by projecting to $\Sigma$ isotopies in $\mathbb{H}^2$ 
for chosen components of the inverse images by $p$ of $\Lambda_\Gamma$ and $E$. 
The obtained graph is isotopic (hence ambient isotopic) to $\Lambda_\Gamma$. Now, 
a vertex $v$ of $\Gamma^d$ which is neither a 
puncture nor a vertex of $\Lambda_\Gamma$ 
lies in a connected component $e$ of $\Sigma \setminus \Gamma$ whose frontier contains 
only non-puncture vertices of $\Gamma$, so according to Lemma~\ref{puncture-graph_lem} it 
lies in the same connected component of $\Sigma \setminus \Lambda_\Gamma$ as $\Gamma$'s 
non-puncture vertices. 
We slide $v$ in $\Sigma \setminus \Lambda_\Gamma$ onto 
its corresponding point in $V$, as follows.   
There is a complementary region $E$ to which $e$ is associated, and  
a connected component $\tilde{E}$ over $E$ is delimited by lifts of curves of 
$\cal T$ (finitely many if $E$ is puncture free). We call $\tilde{F}_E$ the set of lifts of 
frontier paths for $E$ with same endpoints. 
The latter paths determine disks in $\mathbb{H}^2$ over the components of 
$\Sigma \setminus \Gamma$ 
associated to $E$, including $e$, and we take one, $\tilde{e}$, over $e$ 
(if $E$ is puncture free $\tilde{e}$ is unique). 
\begin{figure}
 \centering
 \includegraphics[width=8cm]{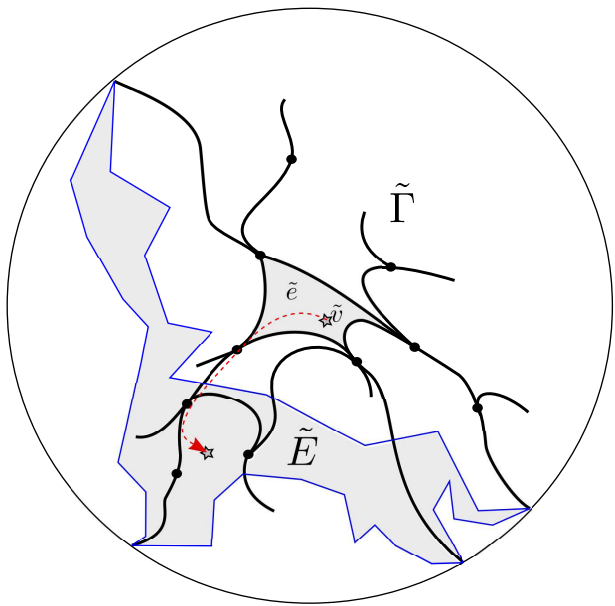}
 \caption{Here $\tilde{e} \cap \tilde{E} = \emptyset$. Four paths in $\tilde{F}_E$ 
determine $\tilde{e}$; a subpath of one of them is projected to $\Sigma$ in order 
to slide $v$ into $E$.}
 \label{Fig.11}
\end{figure}
\noindent
This disk $\tilde{e}$ is 
relatively compact because $e$ is, so the closure of $\tilde{e}$ is a compact disk. 
If $e \cap E \neq \emptyset$ we slide $v$ onto its 
corresponding point in $E$ without exiting $e$. 
If $e \cap E = \emptyset$, then $\tilde{e} \cap \tilde{E} = \emptyset$ 
too, but there is at least one path in $\tilde{F}_E$ among those determining $\tilde{e}$ 
which meets $\tilde{E}$; see Figure~\ref{Fig.11}. 
We use the image by $p$ of a subpath of it to slide $v$ onto 
its corresponding point in $E$; note that since $e$ and $E$ are open sets in $\Sigma$ and a 
frontier path is homotopic to a simple curve (actually a leaf of $\cal T$), the sliding can 
be made along an arc. Any edge of $\Gamma^d$ met during the sliding is dragged 
along. In this way we obtain from $\Gamma^d$ an embedded isotopic graph $\Delta$ with $V$ as 
set of vertices and such that $\Lambda_\Delta \subset (\Sigma \setminus {\cal T})$, as 
wished. 
\EPR

The following lemma gives the desired result for a lamination carried by a filling 
train track. 

\begin{Lem}
Let $\cal T$ be a topological lamination in $\Sigma$ weakly carried by a connected 
train track $\Gamma$ which fills $\Sigma$. Then there is an ambient isotopy of $\Sigma$ 
sending $\cal T$ to a topological lamination contained in a regular neighborhood of 
$\Gamma$, in such a way that $\cal T$ is carried by $\Gamma$ inside this neighborhood 
(i.e. the homotopies deforming $\cal T$'s curves to paths on $\Gamma$ are homotopies 
in the neighborhood). 
\label{Engulfing_fill_lem}
\end{Lem}
\PR
We start with a graph $\Delta$ isotopic to $\Gamma^d$ as constructed in the proof of 
Lemma~\ref{Engulfing_prepa_lem}. 
Then in each complementary region $E$, for each vertex $v$ of $\Delta$ in $E$ 
which is not a vertex of $\Lambda_\Delta$, we choose a PL closed neighborhood of $v$, 
and we choose also a regular neighborhood of $\Lambda_\Delta$ in 
$\Sigma \setminus {\cal T}$. We denote by $N(V)$ the union of all 
these neighborhoods. 
Since $\Delta$ is isotopic to $\Gamma^d$, by Remark~\ref{Fill_dual_reflexive_rem} the 
dual graph $\Delta^d$ is istopic to $\Gamma$, and we can make $\Delta^d$'s 
edges $C^1$, in such a way that the tangent space at each vertex is well-defined and 
compatible with the one at the corresponding vertex of $\Gamma$. 
By Corollary~\ref{Engulfing_fill_cor} the complement $N'$ of the interior of $N(V)$ is a 
regular neighborhood of $\Delta^d$, such that, by construction $\cal T$ is entirely 
contained in the interior of $N'$. By construction too, all homotopies from $\cal T$'s 
curves to paths on $\Delta^d$ can be performed in $N'$. An ambient isotopy taking 
$\Delta^d$ to $\Gamma$ takes $\cal T$ to a topological lamination entirely contained in a 
neighborhood of $\Gamma$. 
\EPR

To settle the general case we need some more technical results. 

\begin{Lem}
Let $\Sigma_{\gamma}$ be the covering space of $\Sigma$ corresponding to the 
cyclic group generated by the homotopy class of an essential simple closed curve 
$\gamma$. Then for every complementary region of $\cal T$, each component 
$\hat{E}$ over $E$ is either a disk or an annulus.
\label{avoid3_lem}
\end{Lem}
\PR
As in Lemma~\ref{avoid1_lem}'s proof,  
the fundamental group of $\Sigma$ is considered at a point $x \in \gamma$, so that 
the covering space $\Sigma_{\gamma}$ corresponding to the cyclic group generated 
by the homotopy class of $\gamma$ is an annulus, and   
the homotopy class of the lift $\hat{\gamma}$ of $\gamma$ to $\Sigma_{\gamma}$ at 
a point $\hat{x}$ over $x$ generates the fundamental group of $\Sigma_{\gamma}$. 
Now let $E$ be a complementary region of $\cal T$ and $\hat{E}$ a component over 
$E$ in $\Sigma_{\gamma}$. The inclusion of $\hat{E}$ in $\Sigma_{\gamma}$ 
induces an injective morphism on fundamental groups, for if not, there should be 
a simple closed curve in $\hat{E}$ bounding a disk in $\Sigma_{\gamma}$ but not 
in $\hat{E}$. But such a disk would contain at least a lift of one of $E$'s 
frontier leaves, which is impossible, so $\hat{E}$ is a disk or an annulus. 
\EPR

\begin{Lem}
Let $\cal T$ be a topological lamination, and $\gamma_c$ an essential PL simple 
closed curve in $\Sigma$ not homotopic to any leaf of $\cal T$. Suppose that 
every leaf of $\cal T$ can be moved away from $\gamma_c$. Then there is a PL 
simple closed curve in $\Sigma \setminus {\cal T}$ homotopic to $\gamma_c$. 
\label{avoid4_lem}
\end{Lem}
\PR
Thanks to Convention~\ref{ClosedLeaf_conv} and Proposition~\ref{LeafByLeaf_prop}, 
there is a geodesic lamination $\cal L$ curve-by-curve isotopic to $\cal T$. Let 
$\gamma_g$ be the geodesic homotopic to $\gamma_c$. Since $\gamma_g$ 
does not belong to $\cal L$ any intersection between $\cal L$ and 
$\gamma_g$ would have to be a transverse one, but if $\cal L$ intersected 
$\gamma_g$, no homotopy on those of $\cal L$'s leaves hitting $\gamma_g$ 
transversally could give curves missing $\gamma_g$, hence $\gamma_c$. Thus 
$\gamma_g$ lies in some complementary region $D$ of $\cal L$, and we denote by 
$E$ the corresponding complementary region of $\cal T$. 
As in Lemma~\ref{avoid1_lem}, $\Sigma_{\gamma_g}$ denotes the covering space of 
$\Sigma$ corresponding to the homotopy class of $\gamma_g$ at a point $x$ chosen 
as base point for the fundamental group. A point $\hat{x}$ is chosen 
over $x$ and $\hat{\gamma_g}$ denotes the lift at $\hat{x}$ of $\gamma_g$.
In $\Sigma_{\gamma_g}$, the component $\hat{D}$ over $D$ is 
either a disk or an annulus, by Lemma~\ref{avoid3_lem}, hence it is an annulus 
because it contains $\hat{\gamma_g}$. 
We take also a component $\tilde{D}$ over $D$ in $\mathbb{H}^2$, and 
denote by $\tilde{\gamma_g}$ a curve over $\hat{\gamma_g}$ in $\tilde{D}$.
We consider 
the component $\tilde{E}$ of $p^{-1}(E)$ in $\mathbb{H}^2$ with the same frontier 
component endpoints on the circle at infinity as $\tilde{D}$, and $\hat{E}$ its 
projection to $\Sigma_{\gamma_g}$. Like $\hat{D}$, $\hat{E}$ is either a disk or 
an annulus.
If the geodesic $\tilde{\gamma_g}$ shares one endpoint with another 
leaf of $\cal L$, both are necessarily frontier leaves of $\cal L$, 
but this contradicts the hypothesis, 
so no lift of a leaf of $\cal L$, hence of $\cal T$, shares an endpoint with 
$\tilde{\gamma_g}$. Now, $E$ has finitely many frontier leaves, and by uniform 
continuity the lift to $\mathbb{H}^2$ of each frontier component of a complementary 
region of $\cal T$, is contained in the $k$-neighborhood of the lift having the same 
endpoints of a frontier leaf of $D$, for some $k > 0$. Let $K$ be the maximum of 
all the $k$ over the frontier leaves of $E$. 
Since each of the frontier leaves of $E$ has its lifts in $\mathbb{H}^2$ at distance 
no more than $K$ from the lift of the frontier component of $D$ sharing the same endpoints, 
and no lift of a leaf of $\cal T$ or $\cal L$ shares an endpoint with $\tilde{\gamma_g}$,
only finitely many of $\tilde{E}$'s frontier components hit $\tilde{\gamma_g}$. 
These frontier components are ambient isotopic in $\mathbb{H}^2$ to the corresponding 
frontier components of $\tilde{D}$, and these ambient isotopies bring $\tilde{\gamma_g}$ 
to a curve $\tilde{\gamma}$ sharing its endpoints with $\tilde{\gamma_g}$ 
and contained in $\tilde{E}$. Note that these ambient isotopies can be undergone in 
$K$-neighborhoods of the involved frontier components of $\tilde{D}$, so that 
if we take in $\mathbb{H}^2$ a (Euclidean) disk centered at $0$ and of (Euclidean) 
radius $r < 1$ sufficiently close to~1, no point of $\tilde{\gamma_g}$ outside that 
disk is contained in any of the above $K$-neighborhoods. In other words the ambient 
isotopies fix the two subarcs of $\tilde{\gamma_g}$ outside the disk. 
Thus, there is a compact subarc 
$a$ of $\tilde{\gamma_g}$ in $\mathbb{H}^2$ such that the above ambient isotopies 
induce a homotopy rel endpoints from $a$ to another arc $a'$ in $\mathbb{H}^2$, and by 
compactness the homotopy is uniformly continuous, so $a'$ 
can be joined to $\tilde{\gamma_g}$'s endpoints by the two subarcs of $\tilde{\gamma_g}$ 
outside the disk,   
to give a curve $\gamma' \subset \tilde{E}$ homotopic to $\tilde{\gamma_g}$. 
Then $\gamma'$ projects to $\hat{E} \subset \Sigma_{\gamma_g}$ as a curve $\hat{\gamma'}$ 
(not necessarily simple or PL, but) homotopic to $\hat{\gamma_g}$, hence  
$\hat{E}$ cannot be a disk, so by Lemma~\ref{avoid3_lem} it is an annulus . We can then 
end the proof as in Lemma~\ref{avoid1_lem}, to obtain a PL simple closed curve in $E$ 
homotopic to $\gamma_g$, hence $\gamma_c$. 
\EPR

\noindent
We have the following corollary: 
\begin{Cor}
Let $\cal T$ be a topological lamination in $\Sigma$ weakly carried by a connected 
train track $\Gamma$ which does not fill $\Sigma$, and let $N$ be a regular 
neighborhood of $\Gamma$. Then there is an ambient isotopy bringing $\cal T$ onto 
a topological lamination none of whose leaves intersects a component of 
$\Sigma \setminus N$ which is neither a disk nor a once punctured disk. 
\label{Avoiding_lam_cor}
\end{Cor}
\PR
Let $\Sigma_0$ be a component of $\Sigma \setminus N$ which is neither a disk nor a 
once punctured disk, and $\gamma_c$ a compact frontier component of $\Sigma_0$. Then 
$\gamma_c$ is essential, and the leaves of $\cal T$ can be moved off $\gamma_c$ 
since they are homotopic to paths on $\Gamma \subset N$. If $\gamma_c$ is homotopic 
to a (necessarily compact) leaf of $\cal T$, an ambient isotopy given 
by~\cite[Th. 2.1]{Eps66} brings the leaf onto $\gamma_c$. Now $\gamma_c$   
being compact, it is a frontier leaf, and it is isolated from at least one 
side because it lies on $\partial N$. We use then a collar neighborhood of 
$\gamma_c$ in the complementary region it is a frontier component of,  
to push $\gamma_c$ into $N$'s interior by an ambient isotopy. 
If $\gamma_c$ is not homotopic to a leaf of $\cal T$, we use  
Lemma~\ref{avoid4_lem} to obtain a simple closed curve 
$\delta$ in some complementary region of $\cal T$, and again,  
an ambient isotopy from $\delta$ to $\gamma_c$ brings $\cal T$ to 
a lamination avoiding $\gamma_c$.
We perform the same operation for each of $\Sigma_0$'s compact frontier components, 
taking care that the isotopies used leave fixed the components already missing 
$\cal T$, until the lamination obtained misses all the compact frontier components of 
$\Sigma_0$. A non-compact frontier component, if any, joins two punctures, and 
is easily seen to be homotopic, hence ambient isotopic to a simple curve \cite{Eps66} in 
$\Sigma \setminus \cal T$ joining the same two punctures (obvious on the universal 
covering), so here too we can move $\cal T$ off such a component by an ambient isotopy. 
Now, the leaves are connected and meet $\Sigma \setminus \Sigma_0$, so since they miss 
its frontier they miss $\Sigma_0$ too. 
We repeat this process with all the components of $\Sigma \setminus N$ as with $\Sigma_0$. 
\EPR

We can now prove the general result about deforming laminations into given regular 
neighborhoods of their carrier graphs.

\begin{Prop}
Let $\cal T$ be a topological lamination in $\Sigma$ weakly carried by a train track 
$\Gamma$. Then there is an ambient isotopy of $\Sigma$ taking $\cal T$ 
to a topological lamination contained in a regular neighborhood of $\Gamma$, in such a 
way that $\cal T$ is carried by $\Gamma$ inside this neighborhood (i.e. the homotopies 
deforming $\cal T$'s curves to paths on $\Gamma$ are homotopies in the neighborhood).
\label{Engulfing_prop}
\end{Prop}
\PR
First, we suppose $\Gamma$ is connected; let $N$ be a regular neighborhood of it. 
We already treated in Lemma~\ref{Engulfing_fill_lem} the case where $\Gamma$ fills 
$\Sigma$, so we now treat the case where some connected components of 
$\Sigma \setminus N$ are neither disks nor once punctured disks, and we apply 
Corollary~\ref{Avoiding_lam_cor} to move $\cal T$ away from all these components.  

In each component of $\Sigma \setminus N$ we choose a point (which is the puncture if the 
component is a once punctured disk), and we can define a graph $\Gamma^d$ dual to 
$\Gamma$ as in the case when $\Gamma$ is filling. However, while in the filling case the 
dual graph is unique up to isotopy in $\Sigma$, this does not hold now, because 
of the components of $\Gamma \setminus N$ which are neither disks nor once punctured disks, 
namely because these components have non cyclic fundamental groups, so a given pair of points 
and/or punctures in the components may be joined to each other by non-homotopic paths. 
Still, a graph which satisfies the two conditions of being dual to $\Gamma^d$, and of 
missing the closure of these components of $\Gamma \setminus N$, is isotopic to $\Gamma$ 
(as a graph). As in the proof of Lemma~\ref{Engulfing_fill_lem}, in the disk 
and once punctured disk components of $\Sigma \setminus N$ we take disjoint PL closed disk 
or once punctured disk neighborhoods of the points chosen above. Then, the complement of 
the union of their interiors together with the other components of $\Sigma \setminus N$ is 
a regular neighborhood $N'$ of a graph $\Delta^d$ dual to $\Gamma^d$. 
Since $\cal T$ misses 
the components of $\Sigma \setminus \Gamma$ which are neither disks nor 
once punctured disks, it is entirely contained in $N'$, and carried by $\Delta^d$,  
as in the proof Lemma~\ref{Engulfing_fill_lem}.

Next, when $\Gamma$ is the union of (finitely many) connected components, each carrying 
one sub-lamination of $\cal T$, two complementary regions of these sub-laminations are 
either nested or disjoint. In particular, we can drag each connected component of $\Gamma$ 
into the innermost complementary region containing the sub-lamination it carries, without 
creating intersections. Note also that lifts of boundary leaves to $\mathbb{H}^2$ do not 
accumulate at each other, except may be on the circle at infinity. Hence we can 
apply the above procedure to each connected component of $\Gamma$ in the 
complementary region it lies in, beginning with the innermost ones. 
\EPR

\section{Appendix}

In this appendix we prove Lemma~\ref{ContinuousToPL_lem}.
\begin{ContToPL_lem}
For any simple infinite curve in a surface $\Sigma$, there is an 
isotopy making it a PL curve. 
If the curve is a one-way infinite curve ending up at a puncture 
in a leaf of a topological lamination, the isotopy can be realized 
by an ambient isotopy avoiding the rest of the lamination. 
\end{ContToPL_lem}
\PR
\noindent
We use the following definition:
A {\smbf square-like} neighborhood of $x \in \Sigma$ is the image $d$ of an 
embedding of the square $[-1, 1] \times [-1, 1]$ into $\Sigma$, such that 
$(0, 0)$ is sent to $x$; we call $x$ the {\smbf center} of $d$ and the 
image of $[-1,1] \times \{0\}$ its {\smbf axis}. 

The first claim we prove is an existence result for square-like neighborhoods: 

\smallskip
\noindent
\underline{Claim 0}: Let $\gamma$ be an infinite curve and $x$ a point on it. 
Then there is a square-like PL neighborhood $d$ with center $x$ and axis a subarc 
of $\gamma$, moreover there is a PL arc in $d$ meeting $\partial d$ only at the 
axis endpoints.
\smallskip

We can suppose $x = \gamma(0)$, and we take a PL disk neighborhood $d$ of $x$ 
in $\Sigma$; then $\gamma^{-1}(int(d))$ is an open set in $\mathbb{R}$ which can be written 
as a (maybe infinite) disjoint union of open intervals, where one of them, J, 
contains~0. Then the two endpoints $y$ and $y'$ of $\gamma(\overline{J})$ lie on 
$\partial d$, and $d$ can be considered as a square-like neighborhood with center 
$x$ and axis $\gamma(\overline{J})$. Since the latter set is compact, we can cover 
it with a finite set of PL disks, and these disks allow us to find a PL arc as we 
want joining $y$ and $y'$ in $d$. This proves Claim~0. 

\smallskip
\noindent
\underline{Claim 1}: Let $d$ be a square-like neighborhood with center $x$ and axis 
$\delta$. Assume there is another arc $\delta'$ in $d$ joining $y = \delta(-1)$ 
to $y' = \delta(1)$, such that $\delta'(0) = \delta(0) = x$, and meeting $\partial d$ 
only at $y$ and $y'$.
Then there is an ambient isotopy from $\delta$ to $\delta'$, and this isotopy can 
be taken in such a way that it pointwise fixes $x, \partial d$, and the complement 
set of $d$. 
\smallskip

\begin{figure}
 \centering
 \includegraphics[width=7.5cm]{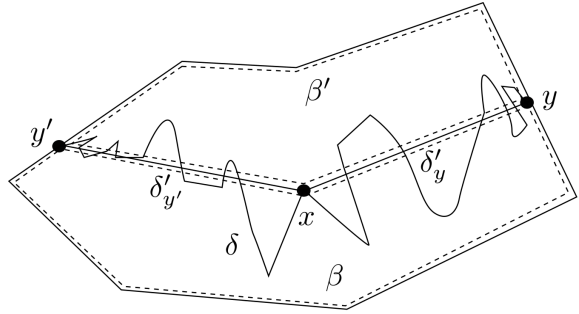}
 \caption{The square-like neighborhood $d$, with 
$(\beta \cup (\delta'_y \cup \delta'_{y'}), \beta' \cup (\delta'_y \cup \delta'_{y'}))$
represented in dashed lines}
 \label{Fig.12}
\end{figure}

\noindent
The proof is an imitation of the first part of the proof of~\cite[Th. A1]{Eps66}.
We denote by $\delta'_y$ and $\delta'_{y'}$ the subarcs of $\delta'$ joining $x$ to 
$y$ and $y'$ respectively. The arcs $\delta'_y$ and $\delta'_{y'}$ together with the 
two arcs joining $y$ and $y'$ in $\partial d$, $\beta$ and $\beta'$, form two pairs 
of circles in $d$, $(\beta \cup \delta, \beta' \cup \delta)$ and 
$(\beta \cup (\delta'_y \cup \delta'_{y'}), \beta' \cup (\delta'_y \cup \delta'_{y'}))$. 
(See Figure~\ref{Fig.12}) 

Now we apply Schoenflies Theorem to conclude that the first components of these two 
pairs bound each one disk in $d$, the two disks being homeomorphic, and the same is 
true for the second components.
These two homeomorphisms fix pointwise $\beta$, $\beta'$ and $x$, and coincide
on $\delta$, so they fit together along $\delta$ to give a homeomorphism $\varphi$ of 
$d$ onto itself, pointwise fixing $\partial d \cup \{x\}$, so we can apply to it 
Alexander's Theorem (also called the {\it Alexander Trick}) to conclude that $\varphi$ is 
isotopic to the identity map of $d$ rel $\partial d \cup \{x\}$, and it sends $\delta$ 
to $\delta'$. This isotopy extends to an ambient isotopy $H$ of $\Sigma$ leaving 
fixed all the points outside of $d$. This proves Claim~1. 
\smallskip

\smallskip
\noindent
\underline{Claim 2}: Let $\delta: J = [r, r'] \rightarrow \Sigma$ be an arc which 
is PL near $r$ and $r'$.
Let $D$ and $D'$ be PL disk neighborhoods of $x = \delta(r)$ and $x' = \delta(r')$ 
respectively, such that $\delta \cap D$ and $\delta \cap D'$ are two PL arcs. 
Then there is an ambient isotopy fixing $D$ and $D'$ from $\delta$ to a PL arc. 
Moreover, given any fixed neighborhood of $\delta$ and any $\varepsilon > 0$, the 
ambient isotopy can be taken in such a way that it fixes the points outside the 
neighborhood, and does not move points by more than $\varepsilon > 0$ inside. 
If $\delta$ misses a compact set then the isotopy can be taken so as to avoid it too. 

\smallskip
\noindent
According to~\cite[Th. A2]{Eps66} applied to the surface with boundary 
$\Sigma \setminus (int(D) \cup int(D'))$, the subarc $\delta^{a}$ of 
$\delta$ between $\partial D$ and $\partial D'$ can be deformed into a PL arc by 
an ambient isotopy rel endpoints. However it seems difficult to estimate how much 
points in $\Sigma$ are moved by this isotopy, so we are going to use here this 
theorem only to conclude that $\delta^{a}$ has a fibered neighborhood $N$ in 
$\Sigma \setminus (int(D) \cup int(D'))$, which can be written as a union of 
$m$ square-like disk neighborhoods $d_i$ satisfying the following four properties:
\begin{itemize}
				\item[P1:] $d_i$ has diameter no more than $\frac{\varepsilon}{2}$
        \item[P2:] $d_1 \cap D$ and $d_m \cap D'$ consist each of a subarc in 
				$\partial D$ and $\partial D'$ respectively 
        \item[P3:] $d_i \cap d_j \neq \emptyset$ iff $|i-j| = 1$, in which case 
				this set consists of one fiber of $N$ 
        \item[P4:] $d_i \cap \delta^{a}$ is the axis of $d_i$
\end{itemize}

\noindent
Such $d_i$ exist for the following reason: As said above, by~\cite[Th. A2]{Eps66} there 
is an ambient isotopy $\phi$ sending $\delta$ to a PL arc with same endpoints, and such a 
PL arc has a PL product neighborhood $N_0 \subset (\Sigma \setminus (int(D) \cup int(D')))$ 
any of whose vertical slicings into disks satisfies 
properties~P2 to~P4. Now, $\delta$ being compact it admits a finite open covering by disks 
of diameter no more that $\frac{\varepsilon}{2}$, and the inverse image by $\phi$ of 
this covering is also an open covering $\cal C$ of $\phi^{-1}(\delta)$, so we can find a 
smaller product neighborhood $N_1$ of $\phi^{-1}(\delta)$ in $N_0$ such that it is covered 
by $\cal C$, 
and then find a vertical slicing of $N_1$ so that each disk produced by this slicing falls 
in the interior of one of the disks of $\cal C$. The image by $\phi$ of the disks produced 
by such a slicing give a set of disks $d_1, \dots, d_m$ (in general not PL ones, however) 
satisfying also property~P1.

Now, the arc $\delta^a$ meets each arc $d_i \cap d_{i+1}$, $i=1, \dots, m-1$ 
in one point $z_i \in int(d_i \cap d_{i+1})$, so there is 
a PL square-like disk $e_i$ with center $z_i$, with axis a subarc of $d_i \cap d_{i+1}$,
and contained in $int(d_i \cup d_{i+1})$. We put $z_0 = \partial D \cap \delta^a$ 
and $z_m = \partial D' \cap \delta^a$. Note also that we can suppose the diameter of 
$e_i$ is no more than $\frac{\varepsilon}{2}$.
We can now apply Claim~1 to $e_i$, $z_i$, and $d_i \cap d_{i+1}$ seen as a curve, to get an 
ambient isotopy such that the axis of $e_i$ goes to a PL arc containing $z_i$ in its interior. 
Doing this for each $i = 1, \dots, m-1$ transforms the $d_i$ into a new 
set of disks $d'_i, i=1, \dots, m$ satisfying properties P1 to P3, and 
satisfying also~P4 with respect to the arc ${\delta'}^a$ obtained from $\delta^a$ 
after having performed the above $m-1$ ambient isotopies. 

Since $z_0$ and $z_m$ lie in two PL subarcs of $\partial D$ and $\partial D'$ respectively, we 
can find in $d'_1$ and $d'_m$ two small PL disks $e_0$ and $e_m$,  
such that $\partial d'_1 \cap e_0$ 
and $\partial d'_m \cap e_m$ are two PL arcs containing respectively $z_0$ and $z_m$ in 
their interiors. Up to taking PL subdisks of the $e_i$ for $i=1, \dots, m-1$, we can assume
that $d'_i \cap e_i$ and $d'_{i+1} \cap e_i$ have the same property with respect to 
$d'_i$ and $d'_{i+1}$.

\begin{figure}
 \centering
 \includegraphics[width=12cm]{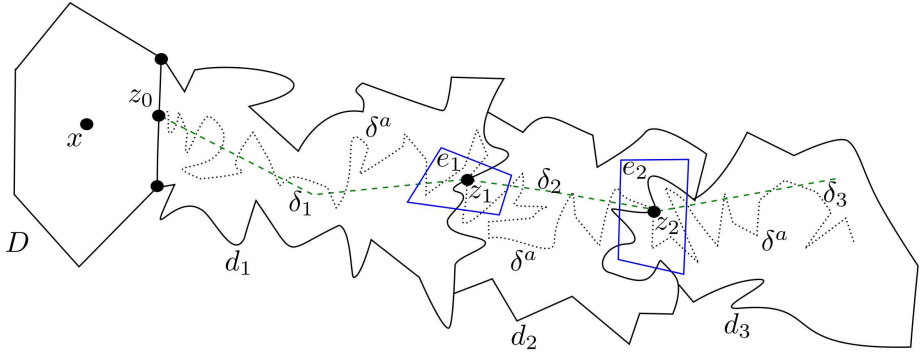}
 \caption{The covering of $\delta^{a}$ by $D$, $D'$ (not represented), and $d_i$}
 \label{Fig.13}
\end{figure}

Now, for each $i = 1, \dots, m$, 
$({\delta'}^a \cap int(d'_i)) \setminus (int(e_{i-1}) \cup int(e_i))$ is a compact set in 
$int(d'_i)$, so it admits a finite covering by PL disks in $int(d'_i)$. These disks 
together with $e_{i-1}$ and $e_i$ form a finite covering of ${\delta'}^a \cap d'_i$, that 
we can use to link $z_{i-1}$ and $z_i$ in $d'_i$ by a simple PL arc $\delta_i$, meeting 
$\partial d'_i$ only along $z_{i-1}$ and $z_i$. (See Figure~\ref{Fig.13})

Applying Claim~1 to $d'_i$ and ${\delta'}^a \cap d'_i$ for each $i$, we get an
ambient isotopy fixing $\partial d'_i$ and sending 
${\delta'}^a \cap d'_i$ to $\delta_i$. This proves Claim~2.

Now we fix an increasing sequence of real numbers 
$\{r_n\}_{n \in \mathbb{Z}}$ or $\{r_n\}_{n \in \mathbb{N}}$ according to whether $\gamma$ is 
two-way infinite or one-way infinite, without accumulation points. 

The strategy is to apply inductively Claims~0 and~1 to $\gamma$, taking the 
indexes of $r_n$ alternatively positive and negative in case $\gamma$ is 
two-way infinite, then Claim~2 to make PL 
the arcs defined on the intervals $[r_j, r_{j+1}]$, $j \in \mathbb{Z}$ or $\mathbb{N}$. 
In the two-way infinite case, we link $k > 0$ and the index of $r_{n_k}$ at 
step~$k$ by the following formula: 
$n_k = (-1)^{k+1} \left\lfloor \frac{k+1}{2} \right\rfloor$. 
We start with a curve $\gamma^{(0)}$, obtained by applying Claims~0 and~1 to $\gamma$ 
and $x^{(0)} = \gamma_i(r_0)$. We denote by $H_0$ the resulting ambient isotopy. 
With the above convention, in any case $n_1 = 1$, and Step~1 consists 
in applying first Claims~0 and~1 to $\gamma^{(0)}$ and $x^{(1)} = \gamma^{(0)}(r_1)$, and 
next Claim~2 to make PL the arc that the resulting curve defines when restricted to 
$[r_0, r_1]$. 
The resulting curve is denoted by $\gamma^{(1)}$. 
When $\gamma$ is two-way infinite (resp. one-way infinite) Step~2 consists in applying 
Claims~0 and~1 to $\gamma^{(1)}$ and $x^{(2)} = \gamma^{(1)}(r_{-1})$ 
(resp. $x^{(2)} = \gamma^{(1)}(r_{2})$), then Claim~2 to make PL the arc that 
the obtained curve defines on $[r_{-1}, r_0]$ (resp. $[r_1, r_2]$), getting 
thus a new curve $\gamma^{(2)}$, and so on. With this notation, the arc of Step~$k \ge 2$ 
to be made PL is the subarc of $\gamma^{(k-1)}$ joining $x^{(k-2)}$ and $x^{(k)}$ 
(resp. $x^{(k-1)}$ and $x^{(k)}$). 
As for the disk sizes, at Step~$k \ge 2$ we require that the disks around the $x^{(k)}$ 
are pairwise disjoint and have diameters no more than $\frac{1}{2^k}$, and also that 
in the isotopy making PL the arcs joining $x^{(k)}$ and $x^{(k-2)}$ (resp. 
$x^{(k-1)}$ and $x^{(k)}$), points are not moved by more than $\frac{1}{2^k}$.
We also ask that the disks at Step~$k$ avoid all the arcs already made PL in the 
preceding steps; this is possible since the union of these arcs is a compact set 
in $\Sigma$.

Applying the above arguments 
inductively produces a sequence of ambient isotopies $H_k$, $k>0$. For $t \in [0, 1]$ fixed, 
the map $F_k(x, t) = H_k(. ,t) \circ \cdots \circ H_{1}(. ,t) \circ H_0(. ,t) (x)$ is a 
homeomorphism, and note that because of our choice of disks, the distance in $\Sigma$ 
between $F_n(x, t)$ and $F_m(x, t)$, $n < m$ is 
bounded above by the sum of the diameters of the disks used between Steps~$n$ and~$m$, that is, 
$\sum_{j=n}^m \frac{1}{2^{j}} < \frac{1}{2^{n-1}}$. This estimation is
independent of $x$ and $t$, so the sequence of functions $F_k(x, t)$ converges uniformly  
and the limit is a map $H$ uniformly continuous on $\Sigma \times [0, 1]$. 
In particular the sequence $\{\gamma^{(k)}\}_{k>0}$ converges, and the 
result is a PL curve $\gamma^{(PL)}$ because after each Step~$k$ all the PL arcs built 
before that step remain fixed in the forthcoming steps. This argument also implies that  
$\gamma^{(PL)}$ is simple since no self-intersection is created at any step. 

We prove now the last assertion. 
We suppose given a topological lamination carried by a train track. First note that a leaf 
ending up at a puncture is isolated from both sides, so given any point on the leaf it is always 
possible to find a neighborhood of the point meeting only the leaf it belongs to. 
For each puncture we fix a nested sequence $\{N_i\}_{i \ge 1}$ of closed once punctured disk 
neighborhoods with empty intersection in $\Sigma$, such that $N_{i+1} \subset int(N_i)$. 
By uniform continuity, given a one-way infinite curve $\gamma$ ending up at 
a puncture in a leaf of the lamination, for each $i$ there is a least value $r_i$ of the 
parameter such that $\gamma([r_i, +\infty)) \subset N_i$. 
Next, we apply Claims~0 to $\gamma$ and the $\gamma(r_i)$, using pairwise disjoint 
square-like neighborhoods, then Claim~1. Since the neighborhoods are disjoint, though they 
are infinite in number the composition of the ambient isotopies given by Claim~1 is an 
ambient isotopy of $\Sigma$. Hence we can assume $\gamma$ is PL around the $\gamma(r_i)$. 
In this situation, for each $i$, $\gamma([0, r_{i+1}])$ is compact, so for each $i$ we 
can find a square-like neighborhood $D_i$ around $\gamma(r_i)$, such that 
$D_i \cap \gamma$ is the PL axis of $D_i$. Hence the conditions are fulfilled to apply 
Claim~2 inductively to the $\delta_i = \gamma([r_i, r_{i+1}])$, with disks $D_i$ and $D_{i+1}$ 
around the endpoints. The difference with the general case, is that we can take the disks 
used in the proof of Claim~2 for $\delta_i^a = \delta_i \setminus (int(D_i \cup D_{i+1}))$ 
disjoint from all the other disks used for $\delta_j^a$, $j < i$. Hence here too, the non-trivial  
effect of the ambient isotopies of Claim~2 takes place in pairwise disjoint 
neighborhoods of the $\delta_i^a$, hence their (infinite) composition is an ambient 
isotopy of $\Sigma$. 
\EPR

{\small \bibliography{global,me}}

\end{document}